\documentclass[reqno,11pt]{amsart}

\usepackage{amsthm,amsmath,amssymb}
\usepackage{mathrsfs,amsfonts,dsfont,functan,extarrows,mathtools}

\usepackage{hyperref}
\usepackage{xcolor}

\newif\ifarxiv
\arxivtrue        

\ifarxiv
  \hypersetup{
      colorlinks=true,
      linkcolor=[rgb]{0.1,0.2,0.6},   
      citecolor=[rgb]{0.1,0.2,0.6},   
      urlcolor=[rgb]{0.1,0.2,0.6},    
      anchorcolor=[rgb]{0.1,0.2,0.6},
      pdfborder={0 0 0},              
      pdftitle={Existence and Geometry of Hermitian Metrics with Constant Second Scalar Curvature},
      pdfauthor={Liangdi Zhang},
      pdfsubject={53C55; 53C21},
      pdfkeywords={Constant second scalar curvature; Compact Hermitian manifolds;  Bismut connection; Chern connection}
  }
\else
  \hypersetup{
      colorlinks=false,               
      pdfborder={0 0 0},              
      pdftitle={Existence and Geometry of Hermitian Metrics with Constant Second Scalar Curvature},
      pdfauthor={Liangdi Zhang}
  }
\fi

\usepackage{marginnote}
\usepackage{xcolor}

\makeatletter
\newcommand{\rmnum}[1]{\romannumeral #1}
\newcommand{\Rmnum}[1]{\expandafter\@slowromancap\romannumeral #1@}
\makeatother

\newtheorem{theorem}{Theorem}[section]
\newtheorem{proposition}[theorem]{Proposition}

\newtheorem{corollary}[theorem]{Corollary}

\newtheorem{lemma}[theorem]{Lemma}

\numberwithin{equation}{section}
\allowdisplaybreaks

\begin{document}

\title[Constant Second Scalar Curvature Metrics]{Existence and Geometry of Hermitian Metrics with Constant Second Scalar Curvature}

\author[L. Zhang]{Liangdi Zhang}
\address[Liangdi Zhang]{\newline Mathematical Science Research Center, Chongqing University of Technology, Chongqing 400054, China}
\email{ldzhang91@163.com}



\begin{abstract}

We study Hermitian metrics with constant second scalar curvature on compact complex manifolds.
We first consider a Yamabe-type problem for the second Bismut scalar curvature within balanced Hermitian conformal classes, and then analyze elliptic equations arising from constant second Chern scalar curvature within a fixed Hermitian conformal class and derive geometric consequences.
Finally, under an Einstein-type condition involving the third and fourth Chern--Ricci curvatures, a pluriclosed Gauduchon Hermitian metric has constant second Chern scalar curvature, which in certain cases further implies the existence of a K\"ahler--Einstein metric.

\vspace*{5pt}



\end{abstract}

\maketitle

\tableofcontents

\section{Introduction}

Constant scalar curvature problems occupy a central position in differential geometry and geometric analysis.
On a compact Riemannian manifold, the classical Yamabe problem asks whether, within a given conformal class, there exists a metric of constant scalar curvature.
Following the foundational contributions of Yamabe~\cite{Yam60}, Trudinger~\cite{Tru68}, Aubin~\cite{Au76}, and Schoen~\cite{Scho84}, it is now known that such metrics always exist.
This celebrated result has motivated a broad range of further investigations.
A related but distinct direction is the non-compact Yamabe problem, which asks whether a smooth complete non-compact Riemannian manifold $(M,g)$ admits a complete metric conformal to $g$ with constant scalar curvature.
In contrast to the compact case, the problem admits a negative answer in general, as shown by counterexamples of Jin~\cite{Jin88}.
Nevertheless, various sufficient conditions guaranteeing existence have been established; see, for instance,
\cite{AMcO88, Gro13, Kim97, Wei19, Hog20}.

The Yamabe problem concerns the prescription of constant scalar curvature within a fixed conformal class.
In a different geometric framework, on compact K\"ahler manifolds one is led to the problem of finding constant scalar curvature K\"ahler (cscK) metrics within a fixed K\"ahler class, which has become a central theme in complex differential geometry.
The cscK equation is a fully nonlinear, fourth-order elliptic equation admitting a variational formulation, whose solutions arise as critical points of the Mabuchi $K$-energy. Moreover, it is deeply connected to algebro-geometric stability through the Yau--Tian--Donaldson conjecture; see, for example, \cite{CDS14a, CDS14b, CDS14c, Tian15, Sze16, BBJ21}. For recent progress on the existence theory of cscK metrics in general K\"ahler classes, we refer to
\cite{ChenCheng21a, ChenCheng21b, Dyrefelt22, ArezzoDVShi23, HanLiu25, Zheng25}.
These developments rely crucially on the K\"ahler condition, which ensures the vanishing of torsion and the equivalence of several canonical connections.

Beyond the K\"ahler category, the situation becomes considerably more subtle.
On general Hermitian manifolds, the presence of torsion leads to several natural scalar curvature notions associated with different canonical connections, notably the Chern and Bismut connections.
In this setting, Angella, Calamai, and Spotti~\cite{ACS17} introduced the Chern--Yamabe problem, which concerns the existence and uniqueness of Hermitian metrics with constant first Chern scalar curvature within a fixed conformal class.
More recently, Barbaro~\cite{Bar23} studied the corresponding Yamabe-type problem associated with the Bismut connection, focusing on the first Bismut scalar curvature and its relation to Calabi--Yau metrics with torsion.

In contrast to the K\"ahler--Einstein case, various Einstein-type conditions associated with different Hermitian connections have been investigated in the Hermitian setting.
The resulting equations are typically non-variational and involve torsion terms, so that the corresponding scalar curvature is not necessarily constant.
For instance, Gauduchon and Ivanov~\cite{GauIva97} studied Einstein--Hermitian and Einstein--Weyl structures in low dimensions, while more recent works considered Einstein equations imposed on the second Chern--Ricci curvature~\cite{Yang25} and on the first Bismut--Ricci curvature together with constant first Bismut scalar curvature~\cite{Ye25}.

Motivated by these developments, we study Hermitian metrics with constant \emph{second} scalar curvature of the Chern and Bismut connections.
Using an analytic approach adapted to Hermitian conformal classes, we establish existence results for the Bismut connection within balanced Hermitian conformal classes, and existence and uniqueness results for the Chern connection under natural curvature assumptions.
We further derive geometric consequences, including restrictions on the Kodaira dimension and on the (anti-)canonical line bundles, as well as Einstein-type conditions for the second Chern scalar curvature that yield constant scalar curvature and criteria related to the existence of K\"ahler--Einstein metrics.

\medskip
To fix notation for later use, let $(M,\omega)$ be a compact Hermitian manifold of complex dimension $n\geq2$. Throughout this paper, all manifolds are assumed to be connected. We denote by $\eta(\omega)$ the \emph{Lee form} associated with $\omega$, which is defined by
\[
d\omega^{n-1}=\eta(\omega)\wedge\omega^{n-1}.
\]

A Hermitian metric $\omega$ is called \emph{balanced} if $d\omega^{n-1}=0$, which is equivalent to $\eta(\omega)=0$, and also to $d^*\omega=0$. It is called a \emph{Gauduchon metric} if
$\partial \bar{\partial} \omega^{n-1} = 0$, which can equivalently be expressed as $\bar{\partial}^* \partial^* \omega = 0$. Alternatively, it is sometimes expressed as $d^*\eta(\omega)=0$ in certain references.
Finally, $\omega$ is \emph{pluriclosed} if
$\partial \bar{\partial} \omega = 0$, and such metrics are also referred to as \emph{strong K\"ahler with torsion} (SKT) in certain references (see, e.g., \cite{FT09}).

It is well known that every Hermitian conformal class $\{\omega\}$ contains a unique Gauduchon metric $\omega_G$, up to scaling (see \cite{Gau77}). Throughout this paper, $\omega_G$ denotes the unique Gauduchon representative in $\{\omega\}$ normalized by $$\mathrm{Vol}(M,\omega_G):=\int_M\frac{\omega_G^n}{n!}=1.$$

A Hermitian conformal class $\{\omega\}$ is said to be \emph{balanced} if its Gauduchon representative $\omega_G$ is balanced.

Let $\Theta^{(i)}(\omega)$, $i=1,2,3,4$, denote the $i$-th Chern--Ricci curvature associated with $\omega$.
We write $S_C^{(1)}(\omega)$ and $S_C^{(2)}(\omega)$ for the first and second Chern scalar curvatures, respectively, and
$S_B^{(1)}(\omega)$ and $S_B^{(2)}(\omega)$ for the first and second Bismut scalar curvatures.
Precise definitions will be given in Section~\ref{sec:2}.

Let $\kappa(M)\in\{-\infty,0,1,\ldots,n\}$ denote the Kodaira dimension of $M$, and let $K_M$ (resp.\ $K_M^{-1}$)
denote the canonical (resp.\ anti-canonical) line bundle of $M$.
Moreover, let $\mathcal{O}_M$ denote the structure sheaf of $M$, that is, the sheaf of holomorphic functions on $M$.

\medskip
We now state the main theorems of this paper.

The first main theorem concerns a Yamabe-type problem for the second Bismut scalar curvature within balanced Hermitian conformal classes.

\begin{theorem}\label{thm1}
Let \((M,\omega)\) be a compact Hermitian manifold of complex dimension
\(n\geq 2\). Assume that the Hermitian conformal class \(\{\omega\}\) is
balanced. Then there exists a Hermitian metric
$\widetilde{\omega}\in\{\omega\}$ with constant second Bismut scalar curvature.
\end{theorem}

It will be seen in Section \ref{sec:4} that if the first Betti number vanishes and the Lee form of the Gauduchon
representative of a given Hermitian conformal class \(\{\omega\}\) is $d$-closed, then \(\{\omega\}\) is balanced. Thus Theorem \ref{thm1} applies.
\begin{corollary}\label{cor4}
Let \((M,\omega)\) be a compact Hermitian manifold of complex dimension
\(n\geq 2\), and let \(\omega_G\in\{\omega\}\) be the Gauduchon
representative. Assume that
\[
b_1(M)=0
\qquad\text{and}\qquad
d\eta(\omega_G)=0.
\]
Then there exists a Hermitian metric $\widetilde{\omega}\in\{\omega\}$ with
constant second Bismut scalar curvature.
\end{corollary}

We turn to the study of the existence and geometry of Hermitian metrics with constant second Chern scalar curvature within a given conformal Hermitian class $\{\omega\}$ on compact complex manifolds. The sign of the total integral of the second Chern scalar curvature of $\omega_G\in\{\omega\}$,
\begin{equation}\label{2TI}
\Gamma_M^{(2)}(\{\omega\}):=\int_MS_C^{(2)}(\omega_G)\frac{\omega_G^n}{n!},
\end{equation}
determines the existence and geometric properties of Hermitian metrics with constant second Chern scalar curvature within the class $\{\omega\}$.
In the special case where the conformal class $\{\omega\}$ is balanced, it also provides information about the first Chern scalar curvature.

\begin{theorem}\label{thm2}
Let $(M,\omega)$ be a compact Hermitian manifold of complex dimension $n\geq 2$ and suppose that $\Gamma_M^{(2)}(\{\omega\}) = 0$. Then there exists a unique (up to scaling) Hermitian metric $\omega_a \in \{\omega\}$ with vanishing second Chern scalar curvature. Moreover, either $\kappa(M)=-\infty$, or $\kappa(M)=0$ and $K_M^{\otimes m}=\mathcal{O}_M$ for some $m\in\mathbb{Z}^+$.

If, in addition, the class $\{\omega\}$ is balanced, then there exists a unique (up to scaling) Hermitian metric $\omega_b \in \{\omega\}$ with vanishing first Chern scalar curvature.
\end{theorem}

\begin{theorem}\label{thm3}
Let $(M,\omega)$ be a compact Hermitian manifold of complex dimension $n\ge 2$, and suppose that
$\Gamma_M^{(2)}(\{\omega\}) < 0$. Then there exists a unique (up to scaling) Hermitian metric \(\omega_c \in \{\omega\}\) with constant second Chern scalar curvature, i.e.,
\[
S_C^{(2)}(\omega_c) = \Gamma_M^{(2)}(\{\omega\}) < 0.
\]

If, in addition, the conformal class \(\{\omega\}\) is balanced, then the anti-canonical line bundle \(K_M^{-1}\) is not pseudo-effective, and there exists a  unique (up to scaling) Hermitian metric \(\omega_d \in \{\omega\}\), such that
\[
S_C^{(1)}(\omega_d) = \Gamma_M^{(2)}(\{\omega\}) < 0.
\]
\end{theorem}

The final main theorem shows that, under natural structural assumptions on a Hermitian metric, an Einstein-type condition associated with the sum of the third and fourth Chern curvatures forces the second Chern scalar curvature to be constant.

\begin{theorem}\label{thm5}
Let $(M,\omega)$ be a compact Hermitian manifold of complex dimension $n\ge 2$, where $\omega$ is both pluriclosed and Gauduchon. Suppose that
\begin{equation}\label{2CE}
\Theta^{(3)}(\omega)+\Theta^{(4)}(\omega)=f\omega
\end{equation}
for some $f\in C^\infty(M,\mathbb{R})$. Then the second Chern scalar curvature of $\omega$ is constant.
\end{theorem}

Recall that in complex dimension two, the pluriclosed condition is equivalent to the Gauduchon condition, while any balanced metric is automatically Gauduchon. As corollaries of Theorem \ref{thm5}, we derive the following rigidity results.

\begin{corollary}\label{cor1}
Let $(M,\omega)$ be a compact Hermitian surface ($n=2$). Suppose that $\omega$ is Gauduchon and satisfies \eqref{2CE}
for some non-positive $f\in C^\infty(M,\mathbb{R})$. Then exactly one of the following holds:
\begin{enumerate}
\item $S_C^{(2)}(\omega)\equiv0$, or
\item $\omega$ is K\"ahler--Einstein with negative scalar curvature.
\end{enumerate}
\end{corollary}

\begin{corollary}\label{cor2}
Let $(M,\omega)$ be a compact Hermitian manifold of complex dimension $n\ge 3$. Suppose  that $\omega$ is pluriclosed and balanced, and satisfies \eqref{2CE}
for some $f\in C^\infty(M,\mathbb{R})$. Then exactly one of the following holds:
\begin{enumerate}
\item $S_C^{(2)}(\omega)\equiv 0$, or
\item $\omega$ is K\"ahler--Einstein with non-zero scalar curvature.
\end{enumerate}
\end{corollary}

As a consequence of Theorems \ref{thm2} and \ref{thm3}, we obtain the following existence and uniqueness result on constant second Chern scalar curvatures.
\begin{corollary}\label{cor3}
Let $(M,\omega)$ be a compact Hermitian manifold of complex dimension $n\ge 2$. If $\omega$ satisfies \eqref{2CE}
for some non-positive $f\in C^\infty(M,\mathbb{R})$, then there exists a unique (up to scaling) $\omega_e\in\{\omega\}$ with constant second Chern scalar curvature, i.e.,
\[
S_C^{(2)}(\omega_e)=\Gamma_M^{(2)}(\{\omega\})\leq0.
\]
\end{corollary}

The paper is organized as follows.
Section~\ref{sec:2} fixes notation and recalls some useful formulas.
Section~\ref{sec:3} is devoted to conformal transformation formulas for Gauduchon metrics and criteria for a Hermitian metric to be K\"ahler.
In Section~\ref{sec:4}, we study a Yamabe-type problem for the Bismut connection, which leads to Theorem~\ref{thm1} and Corollary \ref{cor4}.
Section~\ref{sec:5} addresses the problem of constant second Chern scalar curvature within a fixed conformal Hermitian class, resulting in Theorems~\ref{thm2} and \ref{thm3}.
Section~\ref{sec:6} considers a natural Einstein-type condition for the second Chern scalar curvature, from which Theorem~\ref{thm5} and Corollaries \ref{cor1} to \ref{cor3} are derived. Finally, we present examples of non-K\"ahler metrics with constant Chern scalar curvatures in Section~\ref{sec:7}.
\section{Preliminaries}
\label{sec:2}

Let $(M,\omega)$ be a compact Hermitian manifold of complex dimension $n$, with fundamental $(1,1)$-form
\[
\omega = \sqrt{-1}\, h_{i\bar{j}}\, dz^i \wedge d\bar{z}^j.
\]

Let $g$ denote the underlying Riemannian metric and $J$ the complex structure on $M$.
Then $(M,g,J)$ satisfies
\[
g(X,Y) = g(JX,JY), \quad \omega(X,Y) = g(JX,Y),
\]
for any real vector fields $X,Y \in \mathrm{Sec}(M,T_{\mathbb{R}}M)$.
Moreover, for complex vector fields $W,Z \in \mathrm{Sec}(M,T_{\mathbb{C}}M)$,
\[
g(W,Z) = h(W,Z),
\]
where the complexified tangent bundle is
\[
T_{\mathbb{C}}M = T_{\mathbb{R}}M \otimes \mathbb{C} = T^{1,0}M \oplus T^{0,1}M.
\]

The \emph{Chern connection} ${^C}\nabla$ introduced by Chern \cite{Chern46} is the unique connection compatible with both the Hermitian metric and the holomorphic structure, whose $(0,1)$-part equals $\bar{\partial}$.

Given the Hermitian holomorphic tangent bundle $(T^{1,0}M,h)$, its torsion tensor $T$
and curvature tensor $\Theta$ in local holomorphic coordinates are defined by
\begin{equation}
T_{ij}^k(\omega) = h^{k\bar{l}}\Big(\frac{\partial h_{j\bar{l}}}{\partial z^i} - \frac{\partial h_{i\bar{l}}}{\partial z^j}\Big),
\end{equation}
and
\begin{equation}
\Theta_{i\bar{j}k\bar{l}}(\omega) = -\frac{\partial^2 h_{k\bar{l}}}{\partial z^i \partial \bar{z}^j}
+ h^{p\bar{q}} \frac{\partial h_{p\bar{l}}}{\partial \bar{z}^j} \frac{\partial h_{k\bar{q}}}{\partial z^i},
\end{equation}
respectively.

For a holomorphic vector field $V = V^i \frac{\partial}{\partial z^i} \in \mathrm{Sec}(T^{1,0}M)$,
we define the (1,1)-form
\begin{equation}
T(V) := h_{k\bar{j}} T_{pi}^k(\omega) V^p dz^i \wedge d\bar{z}^j.
\end{equation}

The \emph{Bismut connection} ${^B}\nabla$ \cite{Bismut89, Strominger86}, also called the Strominger or Strominger-Bismut connection,
is the unique Hermitian connection whose torsion ${^B}T$ is totally skew-symmetric. That is,
\[
{^B}\nabla g = 0, \quad {^B}\nabla J = 0, \quad \text{and} \quad {^B}T \in \mathrm{Sec}(\wedge^3 T^*_{\mathbb{R}} M),
\]
where ${^B}T$, regarded as a 3-form, is given by
\[
{^B}T(X,Y,Z) := g\big({^B}\nabla_X Y - {^B}\nabla_Y X - [X,Y], Z\big),
\quad X,Y,Z \in \mathrm{Sec}(T_{\mathbb{R}}M).
\]

Gauduchon \cite{Gau97} introduced a family of Hermitian connections on $(T^{1,0}M,h)$ given by
\begin{equation}
{^t}\nabla=(1-t){^C}\nabla+t{^B}\nabla,\qquad t\in\mathbb{R},
\end{equation}
where the combination is taken in the affine space of Hermitian connections. We call ${^t}\nabla$ the \emph{Gauduchon connections}.

The real curvature tensor of ${^t}\nabla$ on the underlying Riemannian manifold $(M,g)$ is
\[\mathfrak{R}^{t}(X,Y,Z,W)=g({^{t}}\nabla_X{^{t}}\nabla_YZ-{^{t}}\nabla_Y{^{t}}\nabla_XZ-{^{t}}\nabla_{[X,Y]}Z,W),\]
where $X$, $Y$, $Z$, $W\in \mathrm{Sec}(M,T_{\mathbb{R}}M)$. We denote by $R(X,Y,Z,W)(\omega,t)$, where $X$, $Y$, $Z$, $W\in \mathrm{Sec}(M,T_{\mathbb{C}}M)$, and the $\mathbb{C}$-linear complexified curvature tensor on the complexified tangent bundle $T^{\mathbb{C}}M$ is denoted by
\[
R(X,Y,Z,W)(\omega,t), \quad X,Y,Z,W \in \mathrm{Sec}(M,T_{\mathbb{C}}M),
\]
with components
\[
\begin{aligned}
R_{i\bar{j}k\bar{l}}(\omega,t) &= R\Big(\frac{\partial}{\partial z^i},
                                     \frac{\partial}{\partial \bar{z}^j},
                                     \frac{\partial}{\partial z^k},
                                     \frac{\partial}{\partial \bar{z}^l}\Big)(\omega,t),\\[2mm]
R_{ijk\bar{l}}(\omega,t) &= R\Big(\frac{\partial}{\partial z^i},
                                  \frac{\partial}{\partial z^j},
                                  \frac{\partial}{\partial z^k},
                                  \frac{\partial}{\partial \bar{z}^l}\Big)(\omega,t),\\
&\quad \text{etc.}
\end{aligned}
\]

The first, second, third and fourth Ricci curvatures of $(T^{1,0}M,h)$ associated to the Gauduchon connection ${^t}\nabla$ are defined by
\[Ric^{(1)}(\omega,t)=\sqrt{-1}\, R^{(1)}_{i\bar{j}}(\omega,t)dz^i\wedge d\bar{z}^{j}\quad\text{with}\quad R^{(1)}_{i\bar{j}}(\omega,t)=h^{k\bar{l}}R_{i\bar{j}k\bar{l}}(\omega,t),\]
\[Ric^{(2)}(\omega,t)=\sqrt{-1}\, R^{(2)}_{i\bar{j}}(\omega,t)dz^i\wedge d\bar{z}^{j}\quad\text{with}\quad R^{(2)}_{i\bar{j}}(\omega,t)=h^{k\bar{l}}R_{k\bar{l}i\bar{j}}(\omega,t),\]
\[Ric^{(3)}(\omega,t)=\sqrt{-1}\, R^{(3)}_{i\bar{j}}(\omega,t)dz^i\wedge d\bar{z}^{j}\quad\text{with}\quad R^{(3)}_{i\bar{j}}(\omega,t)=h^{k\bar{l}}R_{i\bar{l}k\bar{j}}(\omega,t),\]
\[Ric^{(4)}(\omega,t)=\sqrt{-1}\, R^{(4)}_{i\bar{j}}(\omega,t)dz^i\wedge d\bar{z}^{j}\quad\text{with}\quad R^{(4)}_{i\bar{j}}(\omega,t)=h^{k\bar{l}}R_{k\bar{j}i\bar{l}}(\omega,t),\]
and the first and second scalar curvatures by
\[S^{(1)}(\omega,t)=h^{i\bar{j}}h^{k\bar{l}}R_{i\bar{j}k\bar{l}}(\omega,t)\quad\text{and}\quad S^{(2)}(\omega,t)=h^{i\bar{l}}h^{k\bar{j}}R_{i\bar{j}k\bar{l}}(\omega,t).\]

It is clear that ${^0}\nabla={^C}\nabla$ and ${^1}\nabla={^B}\nabla$. Moreover,
\[
Ric^{(i)}(\omega,0)=\Theta^{(i)}(\omega),\quad Ric^{(i)}(\omega,1)=Ric^{B(i)}(\omega)\quad\text{for}\quad i\in\{1,2,3,4\},
\]
and
\[
S^{(j)}(\omega,0)=S_C^{(j)}(\omega),\quad S^{(j)}(\omega,1)=S_B^{(j)}(\omega)\quad\text{for}\quad j\in\{1,2\}.
\]

For any tensors (or forms) $\alpha$ and $\beta$ of the same bidegree, let $\langle \alpha,\beta \rangle_\omega$ denote their pointwise inner product with respect to the metric $\omega$, and set$|\alpha|^2_\omega=\langle\alpha,\alpha\rangle_\omega$.

Define the global $L^2$-inner product associated with $\omega$ by
\begin{equation}
(\alpha,\beta)_\omega:=\int_M\langle\alpha,\beta\rangle\frac{\omega^n}{n!},\qquad\|\alpha\|^2_\omega:=(\alpha,\alpha)_\omega.
\end{equation}

Let $\Delta_d:=dd^*+d^*d$ be the \emph{Laplace--de Rham operator}, and denote the
\emph{complex Laplacian  operator} associated with $\omega$ by
\[
\Delta_{\omega}^{\mathbb{C}}:=\mathrm{tr}_{\omega}\big(\sqrt{-1}\partial\bar{\partial}(\cdot)\big).
\]

The following formula is frequently used.
\begin{lemma}[see e.g. \cite{Gau84}]
Let $(M,\omega)$ be a compact Hermitian manifold with Lee form $\eta(\omega)$. Then for any $u\in C^\infty(M,\mathbb{R})$, one has
\begin{equation}\label{4.5}
-2\Delta_\omega^\mathbb{C}u=\Delta_du+\langle du,\eta(\omega)\rangle_\omega.
\end{equation}
\end{lemma}

By \cite[Corollary~1.8]{WY25}, for Hermitian metrics, the curvatures of the connections ${^t}\nabla$ and ${^C}\nabla$ are related as follows.
\begin{lemma}[\cite{WY25}]\label{lem1.1}
The curvature tensor of the Gauduchon connection ${^t}\nabla$ on $(M,\omega)$ is
\begin{eqnarray}\label{1.1}
R_{i\bar{j}k\bar{l}}(\omega,t)&=&\Theta_{i\bar{j}k\bar{l}}(\omega)+t(\Theta_{i\bar{l}k\bar{j}}(\omega)+\Theta_{k\bar{j}i\bar{l}}(\omega)-2\Theta_{i\bar{j}k\bar{l}}(\omega))\notag\\
&&+t^2(T_{ik}^p(\omega)\overline{T_{jl}^q}(\omega)h_{p\bar{q}}-h^{p\bar{q}}h_{m\bar{l}}h_{k\bar{r}}T_{ip}^m(\omega)\overline{T_{jq}^r}(\omega)).
\end{eqnarray}

The Ricci curvatures are given by
\begin{equation}\label{1.2}
Ric^{(1)}(\omega,t)=\Theta^{(1)}(\omega)-t(\partial\partial^*\omega+\bar{\partial}\bar{\partial}^*\omega),
\end{equation}
\begin{eqnarray}\label{1.2a}
Ric^{(2)}(\omega,t)&=&\Theta^{(1)}(\omega)-(1-2t)\sqrt{-1}\, \Lambda\partial\bar{\partial}\omega-(1-t)(\partial\partial^*\omega+\bar{\partial}\bar{\partial}^*\omega)\notag\\
&&+(1-t)^2\sqrt{-1}\, T\boxdot \overline{T}-t^2\sqrt{-1}\, T\circ \overline{T},
\end{eqnarray}
\begin{eqnarray}\label{1.2b}
Ric^{(3)}(\omega,t)&=&\Theta^{(1)}(\omega)-t\sqrt{-1}\, \Lambda\partial\bar{\partial}\omega-(1-t)\partial\partial^*\omega-t\bar{\partial}\bar{\partial}^*\omega\notag\\
&&+(t-t^2)\sqrt{-1}\, T\boxdot\overline{T}+t^2T((\partial^*\omega)^\#)
\end{eqnarray}
\begin{eqnarray}\label{1.2c}
Ric^{(4)}(\omega,t)&=&\Theta^{(1)}(\omega)-t\sqrt{-1}\, \Lambda\partial\bar{\partial}\omega-(1-t)\bar{\partial}\bar{\partial}^*\omega-t\partial\partial^*\omega\notag\\
&&+(t-t^2)\sqrt{-1}\, T\boxdot\overline{T}+t^2\overline{T((\partial^*\omega)^\#)},
\end{eqnarray}
where $(\partial^*\omega)^\#$ is the dual holomorphic vector field of the (0,1)-form $\partial^*\omega$, and
\[ T\boxdot\overline{T}:=h^{p\bar{q}}h_{k\bar{l}}T_{ip}^k\cdot\overline{T_{jq}^l}dz^i\wedge d\bar{z}^j,\quad T\circ\overline{T}=h^{p\bar{q}}h^{s\bar{m}}h_{k\bar{j}}h_{i\bar{l}}T_{sp}^k\cdot\overline{T_{mq}^l}dz^i\wedge d\bar{z}^j.\]

The scalar curvatures are related by
\begin{equation}\label{1.3}
S^{(1)}(\omega,t)=S_C^{(1)}(\omega)-2t\langle\partial\partial^*\omega,\omega\rangle_\omega,
\end{equation}
\begin{equation}\label{1.3x}
S^{(2)}(\omega,t)=S_C^{(1)}(\omega)-(1-2t)\langle\partial\partial^*\omega,\omega\rangle_\omega-t^2(2|\partial\omega|^2_\omega+|\partial^*\omega|^2_\omega).
\end{equation}
\end{lemma}

The (first) \emph{Gauduchon degree} of $\{\omega\}$ is defined as the total integral of the first Chern scalar curvature of $\omega_G$, namely,
\begin{equation}
\Gamma_M^{(1)}(\{\omega\}):=\int_MS_C^{(1)}(\omega_G)\frac{\omega_G^n}{n!}=\int_M c_1^{BC}(M)\wedge\frac{\omega_G^{n-1}}{(n-1)!}.
\end{equation}
Motivated by this definition, we call the conformal invariant $\Gamma_M^{(2)}(\{\omega\})$ defined in \eqref{2TI} the \emph{second Gauduchon degree}.

Following \cite{ACS20}, a Hermitian metric $\omega$ is called \emph{weak first Chern--Einstein} if there exists a smooth real-valued function $f$ on $M$ such that
\begin{equation}
\Theta^{(1)}(\omega)=f\omega,
\end{equation}
and $\omega$ is called \emph{weak second Chern--Einstein} if
\begin{equation}
\Theta^{(2)}(\omega)=f\omega.
\end{equation}
This notation is also referred to as \emph{weakly Hermitian--Einstein} in \cite{Yang25}.

Note that neither $\Theta^{(3)}(\omega)$ nor $\Theta^{(4)}(\omega)$ is Hermitian symmetric without additional restriction on $\omega$, whereas their sum $\Theta^{(3)}(\omega)+\Theta^{(4)}(\omega)$ is. This observation motivates the following definition: a Hermitian metric $\omega$ is called a \emph{weak second Hermitian--Einstein metric} if there exists a smooth real-valued function $f$ such that \eqref{2CE} holds.

\section{Conformal transformation formulas and scalar curvatures}
\label{sec:3}

In this section, we derive the conformal transformation formulas for the curvatures associated with the Gauduchon connection ${^t}\nabla$ under conformal changes of a Hermitian metric $\omega$.
These formulas play a fundamental role in the analysis of conformal equations arising in the study of constant second Bismut and Chern scalar curvature metrics.

After deriving these formulas, we turn to the comparison between the first and second total Gauduchon scalar curvatures, and show that such a comparison yields criteria for a Hermitian metric to be K\"ahler.

\begin{proposition}\label{prp3.1}
Define $\omega_f=e^f\omega$ for some $f\in C^\infty(M,\mathbb{R})$. The third Ricci curvature of the Gauduchon connection ${^t}\nabla$ on $(M,\omega_f)$ satisfies
\begin{eqnarray}\label{1.4}
Ric^{(3)}(\omega_f,t)&=&Ric^{(3)}(\omega,t)-(1+(n-2)t)\sqrt{-1}\, \partial\bar{\partial}f-t\Delta_\omega^{\mathbb{C}}f\cdot\omega\notag\\
&&-nt^2|\partial f|^2_\omega\cdot\omega+t^2\sqrt{-1}\, \partial f\wedge\bar{\partial}f\notag\\
&&-nt^2\sqrt{-1}\, T((\bar{\partial} f)^\#)-t^2\overline{\sqrt{-1}\, T((\bar{\partial} f)^\#)}\notag\\
&&+t^2\langle\partial^*\omega,\sqrt{-1}\, \bar{\partial}f\rangle_\omega\cdot\omega-t^2\bar{\partial}^*\omega\wedge\bar{\partial}f,
\end{eqnarray}
and the second scalar curvature is
\begin{eqnarray}\label{1.5}
S^{(2)}(\omega_f,t)&=&e^{-f}\big(S^{(2)}(\omega,t)-(1+2(n-1)t)\Delta^{\mathbb{C}}_{\omega} f\notag\\
&&-(n^2-1)t^2|\partial f|^2_\omega+2(n+1)t^2\mathrm{Re}\langle\partial^*\omega,\sqrt{-1}\bar{\partial}f\rangle_\omega\big).
\end{eqnarray}
\end{proposition}

\begin{proof}
It follows from \eqref{1.1} that
\begin{eqnarray}\label{1.6}
&&R_{i\bar{j}k\bar{l}}(\omega_f,t)\notag\\
&=&\Theta_{i\bar{j}k\bar{l}}(\omega_f)+t(\Theta_{i\bar{l}k\bar{j}}(\omega_f)+\Theta_{k\bar{j}i\bar{l}}(\omega_f)-2\Theta_{i\bar{j}k\bar{l}}(\omega_f))\notag\\
&&+t^2(T_{ik}^p(\omega_f)\overline{T_{jl}^q}(\omega_f)(h_f)_{p\bar{q}}-(h_f)^{p\bar{q}}(h_f)_{m\bar{l}}(h_f)_{k\bar{r}}T_{ip}^m(\omega_f)\overline{T_{jq}^r}(\omega_f)),
\end{eqnarray}
where $(h_f)_{i\bar{j}}=e^{f}h_{i\bar{j}}$ and $(h_f)^{i\bar{j}}=e^{-f}h^{i\bar{j}}$.

One can compute that (see e.g. \cite[(2.10)]{YZ25})
\begin{equation}\label{1.7}
\Theta_{i\bar{j}k\bar{l}}(\omega_f)=e^f(\Theta_{i\bar{j}k\bar{l}}(\omega)-h_{k\bar{l}}\frac{\partial^2f}{\partial z^i\partial\bar{z}^j}),
\end{equation}
and then
\begin{eqnarray}\label{1.8}
&&\Theta_{i\bar{l}k\bar{j}}(\omega_f)+\Theta_{k\bar{j}i\bar{l}}(\omega_f)-2\Theta_{i\bar{j}k\bar{l}}(\omega_f)\notag\\
&=&e^f(\Theta_{i\bar{l}k\bar{j}}(\omega)+\Theta_{k\bar{j}i\bar{l}}(\omega)-2\Theta_{i\bar{j}k\bar{l}}(\omega))\notag\\
&&-e^f(h_{k\bar{j}}\frac{\partial^2f}{\partial z^i\partial\bar{z}^l}+h_{i\bar{l}}\frac{\partial^2f}{\partial z^k\partial\bar{z}^j}-2h_{k\bar{l}}\frac{\partial^2f}{\partial z^i\partial\bar{z}^j}).
\end{eqnarray}

Note that
\begin{eqnarray*}
T_{ij}^k(\omega_f)&=&(h_f)^{k\bar{l}}\frac{\partial(h_f)_{j\bar{l}}}{\partial z^i}-(h_f)^{k\bar{l}}\frac{\partial(h_f)_{i\bar{l}}}{\partial z^j}\notag\\
&=&h^{k\bar{l}}\frac{\partial h_{j\bar{l}}}{\partial z^i}+h^{k\bar{l}}h_{j\bar{l}}\frac{\partial f}{\partial z^i}-h^{k\bar{l}}\frac{\partial h_{i\bar{l}}}{\partial z^j}-h^{k\bar{l}}h_{i\bar{l}}\frac{\partial f}{\partial z^j}\notag\\
&=&T_{ij}^k(\omega)+\delta_j^k\frac{\partial f}{\partial z^i}-\delta_i^k\frac{\partial f}{\partial z^j}.
\end{eqnarray*}

Therefore,
\begin{eqnarray*}
&&T_{ik}^p(\omega_f)\overline{T_{jl}^q}(\omega_f)(h_f)_{p\bar{q}}\notag\\
&=&(T_{ik}^p(\omega)+\delta_k^p\frac{\partial f}{\partial z^i}-\delta_i^p\frac{\partial f}{\partial z^k})(\overline{T_{jl}^q}(\omega)+\delta_{\bar{l}}^{\bar{q}}\frac{\partial f}{\partial \bar{z}^j}-\delta_{\bar{j}}^{\bar{q}}\frac{\partial f}{\partial \bar{z}^l})e^fh_{p\bar{q}}\notag\\
&=&e^fT_{ik}^p(\omega)\overline{T_{jl}^q}(\omega)h_{p\bar{q}}+e^f(h_{k\bar{q}}\frac{\partial f}{\partial z^i}-h_{i\bar{q}}\frac{\partial f}{\partial z^k})\overline{T_{jl}^q}(\omega)\notag\\
&&+e^fT_{ik}^p(\omega)(h_{p\bar{l}}\frac{\partial f}{\partial \bar{z}^j}-h_{p\bar{j}}\frac{\partial f}{\partial \bar{z}^l})\notag\\
&&+e^f(h_{k\bar{l}}\frac{\partial f}{\partial z^i}\frac{\partial f}{\partial \bar{z}^j}-h_{i\bar{l}}\frac{\partial f}{\partial z^k}\frac{\partial f}{\partial \bar{z}^j}-h_{k\bar{j}}\frac{\partial f}{\partial z^i}\frac{\partial f}{\partial \bar{z}^l}+h_{i\bar{j}}\frac{\partial f}{\partial z^k}\frac{\partial f}{\partial \bar{z}^l}),
\end{eqnarray*}
and
\begin{eqnarray*}
&&(h_f)^{p\bar{q}}(h_f)_{m\bar{l}}(h_f)_{k\bar{r}}T_{ip}^m(\omega_f)\overline{T_{jq}^r}(\omega_f)\notag\\
&=&e^fh^{p\bar{q}}h_{m\bar{l}}h_{k\bar{r}}(T_{ip}^m(\omega)+\delta_p^m\frac{\partial f}{\partial z^i}-\delta_i^m\frac{\partial f}{\partial z^p})(\overline{T_{jq}^r}(\omega)+\delta_{\bar{q}}^{\bar{r}}\frac{\partial f}{\partial \bar{z}^j}-\delta_{\bar{j}}^{\bar{r}}\frac{\partial f}{\partial \bar{z}^q})\notag\\
&=&e^fh^{p\bar{q}}h_{m\bar{l}}h_{k\bar{r}}T_{ip}^m(\omega)\overline{T_{jq}^r}(\omega)+e^fh_{k\bar{r}}\overline{T_{jl}^r}(\omega)\frac{\partial f}{\partial z^i}-e^fh^{p\bar{q}}h_{i\bar{l}}h_{k\bar{r}}\frac{\partial f}{\partial z^p}\overline{T_{jq}^r}(\omega)\notag\\
&&+e^fh_{m\bar{l}}T_{ik}^m(\omega)\frac{\partial f}{\partial \bar{z}^j}+e^fh_{k\bar{l}}\frac{\partial f}{\partial z^i}\frac{\partial f}{\partial \bar{z}^j}-e^fh_{i\bar{l}}\frac{\partial f}{\partial z^k}\frac{\partial f}{\partial \bar{z}^j}\notag\\
&&-e^fh^{p\bar{q}}h_{m\bar{l}}h_{k\bar{j}}T_{ip}^m(\omega)\frac{\partial f}{\partial \bar{z}^q}-e^fh_{k\bar{j}}\frac{\partial f}{\partial z^i}\frac{\partial f}{\partial \bar{z}^l}+e^fh_{i\bar{l}}h_{k\bar{j}}h^{p\bar{q}}\frac{\partial f}{\partial z^p}\frac{\partial f}{\partial \bar{z}^q}.
\end{eqnarray*}
It follows that
\begin{eqnarray}\label{1.12}
&&T_{ik}^p(\omega_f)\overline{T_{jl}^q}(\omega_f)(h_f)_{p\bar{q}}-(h_f)^{p\bar{q}}(h_f)_{m\bar{l}}(h_f)_{k\bar{r}}T_{ip}^m(\omega_f)\overline{T_{jq}^r}(\omega_f)\notag\\
&=&e^f(T_{ik}^p(\omega)\overline{T_{jl}^q}(\omega)h_{p\bar{q}}-h^{p\bar{q}}h_{m\bar{l}}h_{k\bar{r}}T_{ip}^m(\omega)\overline{T_{jq}^r}(\omega))\notag\\
&&-e^fh_{i\bar{q}}\frac{\partial f}{\partial z^k}\overline{T_{jl}^q}(\omega)-e^fh_{p\bar{j}}\frac{\partial f}{\partial \bar{z}^l}T_{ik}^p(\omega)+e^fh_{i\bar{j}}\frac{\partial f}{\partial z^k}\frac{\partial f}{\partial \bar{z}^l}\notag\\
&&+e^fh^{p\bar{q}}h_{i\bar{l}}h_{k\bar{r}}\frac{\partial f}{\partial z^p}\overline{T_{jq}^r}(\omega)+e^fh^{p\bar{q}}h_{m\bar{l}}h_{k\bar{j}}T_{ip}^m(\omega)\frac{\partial f}{\partial \bar{z}^q}\notag\\
&&-e^fh_{i\bar{l}}h_{k\bar{j}}h^{p\bar{q}}\frac{\partial f}{\partial z^p}\frac{\partial f}{\partial \bar{z}^q}.
\end{eqnarray}

Applying \eqref{1.7}, \eqref{1.8} and \eqref{1.12} to \eqref{1.6}, we get
\begin{eqnarray}\label{1.15}
&&R_{i\bar{l}}^{(3)}(\omega_f,t)\notag\\
&=&(h_f)^{k\bar{j}}R_{i\bar{j}k\bar{l}}(\omega_f,t)\notag\\
&=&h^{k\bar{j}}\Theta_{i\bar{j}k\bar{l}}(\omega)+th^{k\bar{j}}(\Theta_{i\bar{l}k\bar{j}}(\omega)+\Theta_{k\bar{j}i\bar{l}}(\omega)-2\Theta_{i\bar{j}k\bar{l}}(\omega))\notag\\
&&+t^2h^{k\bar{j}}(T_{ik}^p(\omega)\overline{T_{jl}^q}(\omega)h_{p\bar{q}}-h^{p\bar{q}}h_{m\bar{l}}h_{k\bar{r}}T_{ip}^m(\omega)\overline{T_{jq}^r}(\omega))\notag\\
&&-h^{k\bar{j}}h_{k\bar{l}}\frac{\partial^2f}{\partial z^i\partial\bar{z}^j}-th^{k\bar{j}}(h_{k\bar{j}}\frac{\partial^2f}{\partial z^i\partial\bar{z}^l}+h_{i\bar{l}}\frac{\partial^2f}{\partial z^k\partial\bar{z}^j}-2h_{k\bar{l}}\frac{\partial^2f}{\partial z^i\partial\bar{z}^j})\notag\\
&&+t^2h^{k\bar{j}}(-h_{i\bar{q}}\frac{\partial f}{\partial z^k}\overline{T_{jl}^q}(\omega)-h_{p\bar{j}}\frac{\partial f}{\partial \bar{z}^l}T_{ik}^p(\omega)+h_{i\bar{j}}\frac{\partial f}{\partial z^k}\frac{\partial f}{\partial \bar{z}^l})\notag\\
&&+t^2h^{k\bar{j}}(h^{p\bar{q}}h_{i\bar{l}}h_{k\bar{r}}\frac{\partial f}{\partial z^p}\overline{T_{jq}^r}(\omega)+h^{p\bar{q}}h_{m\bar{l}}h_{k\bar{j}}T_{ip}^m(\omega)\frac{\partial f}{\partial \bar{z}^q}-h_{i\bar{l}}h_{k\bar{j}}h^{p\bar{q}}\frac{\partial f}{\partial z^p}\frac{\partial f}{\partial \bar{z}^q})\notag\\
&=&h^{k\bar{l}}R_{i\bar{j}k\bar{l}}(\omega,t)-(1+(n-2)t)\frac{\partial^2f}{\partial z^i\partial\bar{z}^l}-th^{k\bar{j}}\frac{\partial^2f}{\partial z^k\partial\bar{z}^j}h_{i\bar{l}}\notag\\
&&+t^2\big(-h^{k\bar{j}}h_{i\bar{q}}\overline{T_{jl}^q}(\omega)\frac{\partial f}{\partial z^k}-\sum_pT_{ip}^p(\omega)\frac{\partial f}{\partial\bar{z}^l}+\frac{\partial f}{\partial z^i}\frac{\partial f}{\partial\bar{z}^l}\big)\notag\\
&&-t^2(h^{k\bar{j}}h_{i\bar{l}}\sum_{\bar{q}}\overline{T_{jq}^q}(\omega)\frac{\partial f}{\partial z^k}+nh^{k\bar{j}}h_{p\bar{l}}T_{ki}^p(\omega)\frac{\partial f}{\partial\bar{z}^j}+nh^{k\bar{j}}\frac{\partial f}{\partial z^k}\frac{\partial f}{\partial\bar{z}^j}h_{i\bar{l}}\big).
\end{eqnarray}

It is known that (see e.g. \cite[Lemma~3.3]{LY17})
\begin{equation}
\partial^*\omega=-\sqrt{-1}\, \sum_{\bar{q}}\overline{T_{jq}^q}(\omega)d\bar{z}^j\quad\text{and}\quad\bar{\partial}^*\omega=\sqrt{-1}\, \sum_pT_{ip}^p(\omega)dz^i.
\end{equation}
Then we have
\begin{equation}\label{1.16}
\sqrt{-1}\, \sum_pT_{ip}^p(\omega)\frac{\partial f}{\partial\bar{z}^l}dz^i\wedge d\bar{z}^l=\bar{\partial}^*\omega\wedge\bar{\partial}f,
\end{equation}
and
\begin{equation}\label{1.17}
-h^{k\bar{j}}\sum_{\bar{q}}\overline{T_{jq}^q}(\omega)\frac{\partial f}{\partial z^k}=\langle\partial^*\omega,\sqrt{-1}\, \bar{\partial}f\rangle_\omega.
\end{equation}

Applying \eqref{1.16} and \eqref{1.17} to \eqref{1.15}, we get \eqref{1.4}.

Taking trace of \eqref{1.15}, we obtain
\begin{eqnarray}\label{1.13}
&&S^{(2)}(\omega_f,t)\notag\\
&=&(h_f)^{i\bar{l}}Ric^{(3)}_{i\bar{l}}(\omega_f,t)\notag\\
&=&e^{-f}\big(S^{(2)}(\omega,t)-(1+2(n-1)t)\Delta^{\mathbb{C}}_{\omega} f)-t^2(n^2-1)|\partial f|_\omega^2\notag\\
&&-(n+1)t^2(h^{k\bar{j}}\frac{\partial f}{\partial z^k}\sum_{\bar{q}}\overline{T_{jq}^q}(\omega)+h^{k\bar{j}}\frac{\partial f}{\partial\bar{z}^j}\sum_pT_{kp}^p(\omega)\big)\notag\\
&=&e^{-f}\big(S^{(2)}(\omega,t)-(1+2(n-1)t)\Delta^{\mathbb{C}}_{\omega} f\notag\\
&&-(n^2-1)t^2|\partial f|^2_\omega+2(n+1)t^2\mathrm{Re}\langle\partial^*\omega,\sqrt{-1}\bar{\partial}f\rangle_\omega\big).
\end{eqnarray}
This is \eqref{1.5}.
\end{proof}

Taking $t=0$ and $t=1$ in Proposition \ref{prp3.1}, we respectively derive the conformal formulas for the scalar curvatures associated with the Chern connection and the Bismut connection.
\begin{corollary}
Let $(M,\omega)$ be a compact Hermitian manifold of complex dimension $n$, and let $\omega_f=e^f\omega$ for some $f\in C^\infty(M,\mathbb{R})$. Then the following conformal transformation formulas hold.

(\rmnum{1}) The third and fourth Chern--Ricci curvatures satisfy
\begin{equation}\label{4.1}
\Theta^{(3)}(\omega_f)=\Theta^{(3)}(\omega)-\sqrt{-1}\, \partial\bar{\partial}f,
\end{equation}
and
\begin{equation}\label{4.1x}
\Theta^{(4)}(\omega_f)=\Theta^{(4)}(\omega)-\sqrt{-1}\, \partial\bar{\partial}f.
\end{equation}

(\rmnum{2}) The second Chern scalar curvature transforms as
\begin{eqnarray}\label{4.2}
S_C^{(2)}(\omega_f)&=&e^{-f}(S_C^{(2)}(\omega)-\Delta^{\mathbb{C}}_{\omega} f).
\end{eqnarray}

(\rmnum{3}) The third Bismut--Ricci curvature is given by
\begin{eqnarray}\label{4.3}
Ric^{B(3)}(\omega_f)&=&Ric^{B(3)}(\omega)-(n-1)\sqrt{-1}\, \partial\bar{\partial}f-\Delta_\omega^{\mathbb{C}}f\cdot\omega\notag\\
&&-n|\partial f|^2_\omega\cdot\omega+\sqrt{-1}\, \partial f\wedge\bar{\partial}f\notag\\
&&-n\sqrt{-1}\, T((\bar{\partial} f)^\#)-\overline{\sqrt{-1}\, T((\bar{\partial} f)^\#)}\notag\\
&&+\langle\partial^*\omega,\sqrt{-1}\, \bar{\partial}f\rangle_\omega\cdot\omega-\bar{\partial}^*\omega\wedge\bar{\partial}f,
\end{eqnarray}

(\rmnum{4}) The second Bismut scalar curvature satisfies
\begin{eqnarray}\label{4.4}
S_B^{(2)}(\omega_f)&=&e^{-f}\big(S_B^{(2)}(\omega)-(2n-1)\Delta^{\mathbb{C}}_{\omega} f\notag\\
&&-(n^2-1)|\partial f|^2_\omega+2(n+1)\mathrm{Re}\langle\partial^*\omega,\sqrt{-1}\bar{\partial}f\rangle_\omega\big).
\end{eqnarray}
\end{corollary}

We turn to the comparison between the first and second total Gauduchon scalar curvatures.
\begin{proposition}\label{prp7.1}
Let $M$ be a complete complex manifold equipped with a Gauduchon metric $\omega$. Then
\begin{equation}\label{7.1}
S^{(2)}(\omega,t) = S^{(1)}(\omega,t) - (t^2 - 4t + 1)|\bar{\partial}^*\omega|_\omega^2 - 2t^2|\partial\omega|_\omega^2.
\end{equation}

In particular, in complex dimension two, one has
\begin{equation}\label{7.2}
S^{(2)}(\omega,t) = S^{(1)}(\omega,t) - (3t-1)(t-1)|\partial\omega|_\omega^2.
\end{equation}
\end{proposition}

\begin{proof}
Equation \eqref{7.1} follows directly from \eqref{1.3}, \eqref{1.3x}, and \cite[Lemma 4.5]{Yang25b+}.
Note that
\[
|\partial\omega|^2 = |*\partial*\omega|^2 = |\bar{\partial}^*\omega|^2
\]
in complex dimension two, \eqref{7.2} follows.
\end{proof}

\begin{theorem}
Let $(M,\omega)$ be a complete Hermitian manifold of complex dimension $n$ equipped with a Gauduchon metric satisfying
\begin{equation}\label{7.3}
-\infty<\int_MS^{(1)}(\omega,t)\frac{\omega^n}{n!}\leq \int_MS^{(2)}(\omega,t)\frac{\omega^n}{n!}<+\infty.
\end{equation}

\begin{enumerate}
\item If $n\geq3$ and $t=0$, then $\omega$ is balanced and $S^{(1)}(\omega,t) = S^{(2)}(\omega,t)$.
\item If $n\geq3$, $t\neq0$ and $\omega$ is balanced, then $\omega$ is K\"ahler.
\item If $n\geq3$ and $t\in(-\infty,0)\cup(0,2-\sqrt{3}]\cup[2+\sqrt{3},+\infty)$, then $\omega$ is K\"ahler.
\item If $n=2$ and $t\in(-\infty,\frac{1}{3})\cup(1,+\infty)$, then $\omega$ is K\"ahler.
\end{enumerate}
\end{theorem}

\begin{proof}
Applying \eqref{7.1} to case (1), we obtain $\|\bar{\partial}^*\omega\|_\omega^2 = 0$, which implies that $\omega$ is balanced. Consequently, $S^{(1)}(\omega,t) = S^{(2)}(\omega,t)$.

In cases (2) and (3), \eqref{7.1} gives $\|\partial\omega\|_\omega^2 = 0$, and therefore $\omega$ is K\"ahler.

In case (4), \eqref{7.2} similarly implies $\|\partial\omega\|_\omega^2 = 0$, and hence $\omega$ is K\"ahler.
\end{proof}

\section{A Yamabe-type problem for the Bismut connection}
\label{sec:4}

Let \((M,\omega)\) be a compact Hermitian manifold of complex dimension
\(n\geq 2\), and assume that the Hermitian conformal class
\(\{\omega\}\) is balanced. Let \(\omega_G\in\{\omega\}\) be its
Gauduchon representative. By definition, \(\omega_G\) is balanced.
Set
\[
\omega_B:=\omega_G.
\]
Then
\[
\eta(\omega_B)=0,
\qquad
\partial^*\omega_B=0.
\]

In this section,
we prove that there exists $f\in C^\infty(M,\mathbb{R})$ so that $\omega_f:=e^{f}\omega_B$ has constant second Bismut scalar curvature, i.e., $$S_B^{(2)}(\omega_f)=\lambda\quad\text{for some real constant } \lambda.$$
By \eqref{4.4} and $\partial^*\omega_B=0$, this problem reduces to solving the semilinear elliptic equation of
\begin{equation}\label{5.1}
(2n-1)\Delta^{\mathbb{C}}_{\omega_B} f+(n^2-1)|\partial f|^2_{\omega_B}-S_B^{(2)}(\omega_B)=-\lambda e^{f}.
\end{equation}

Let
\[
\varphi:=\exp\{\frac{n^2-1}{2n-1} f\},
\]
\eqref{5.1} is equivalent to
\begin{eqnarray}\label{5.2}
\square_{\omega_B}\varphi:=-\Delta^{\mathbb{C}}_{\omega_B} \varphi+N_1S_B^{(2)}(\omega_B)\varphi=N_1\lambda\varphi^{N_2-1},
\end{eqnarray}
where
$$N_1=\frac{n^2-1}{(2n-1)^2}\quad\text{and}\quad N_2=2+\frac{2n-1}{n^2-1}.$$

For $2< q< \frac{2n}{n-1}$ ($n\geq2$) and real-valued functions $0\not\equiv \varphi\in W^{1,2}(M,\omega_B)$, we consider the functional
\begin{equation}\label{5.3}
Y_q(\varphi):=\int_M\varphi \square_{\omega_B}\varphi\frac{\omega_B^n}{n!}\cdot\big(N_1\int_M|\varphi|^q\frac{\omega_B^n}{n!}\big)^{-\frac{2}{q}},
\end{equation}
and define
\begin{equation}\label{5.4}
\mu_q:=\inf \{Y_q(\varphi)\, |\, 0\not\equiv\varphi\in W^{1,2}(M,\omega_B)\}.
\end{equation}

For the sake of simplicity, we denote the $L^p(M,\omega_B)$-norm of a $L^p$ function on $M$ with respect to the metric $\omega_B$ by $\|\cdot\|_p$ and the $W^{k,p}(M,\omega_B)$-norm by $\|\cdot\|_{k,p}$.
\begin{proposition}\label{prp5.2}
For any $2< q< \frac{2n}{n-1}$ ($n\geq2$),
\begin{equation}\label{5.3x}
Y_q(\varphi)=N_1^{-\frac{2}{q}}\big(\|\partial\varphi\|_2^2+N_1\int_MS_B^{(2)}(\omega_B)\varphi^2\frac{\omega_B^n}{n!}\big)\cdot\|\varphi\|_q^{-2}
\end{equation}
and $\mu_q$ is uniformly finite.
\end{proposition}

\begin{proof}
According to Sobolev embedding theorem for compact manifolds, one has $W^{1,2}(M,\omega_B)$ is continuously embedded in $L^q(M,\omega_B)$. It ensures that the denominator of $Y_q(\varphi)$ makes sense.

Using \eqref{4.5} and $\eta(\omega_B)=0$, we obtain
\begin{eqnarray}\label{5.6}
\int_M\varphi \square_{\omega_B}\varphi\frac{\omega_B^n}{n!}&=&\frac{1}{2}\int_M\varphi\Delta_d\varphi\frac{\omega_B^n}{n!}+N_1\int_MS_B^{(2)}(\omega_B)\varphi^2\frac{\omega_B^n}{n!}\notag\\
&=&\|\partial\varphi\|_2^2+N_1\int_MS_B^{(2)}(\omega_B)\varphi^2\frac{\omega_B^n}{n!},
\end{eqnarray}

\eqref{5.3x} follows by applying \eqref{5.6} to \eqref{5.3}.

Since \(\omega_B=\omega_G\) and \(\operatorname{Vol}(M,\omega_G)=1\), we have
\begin{eqnarray}\label{eq4.7}
\mu_q&\leq&Y_q(1)\notag\\
&=&N_1^{1-\frac{2}{q}}\int_M S_B^{(2)}(\omega_B)\frac{\omega_B^n}{n!}\notag\\
&\leq&\sup_M |S_B^{(2)}(\omega_B)|.
\end{eqnarray}

On the other hand, if \(\inf_M S_B^{(2)}(\omega_B)\geq0\), then \(\mu_q\geq0\) by
\eqref{5.3x}. If \(\inf_M S_B^{(2)}(\omega_B)<0\), then H\"older's inequality gives
\[
\|\varphi\|_2\leq\|\varphi\|_q,
\]
and hence
\[
Y_q(\varphi)
\geq
N_1^{1-\frac{2}{q}}
\inf_M S_B^{(2)}(\omega_B)\cdot
\|\varphi\|_2^2\|\varphi\|_q^{-2}
\geq
\inf_M S_B^{(2)}(\omega_B).
\]
Consequently,
\begin{equation}\label{eq4.8}
\min\left\{\inf_M S_B^{(2)}(\omega_B),0\right\}
\leq\mu_q\leq\sup_M|S_B^{(2)}(\omega_B)|.
\end{equation}
Thus \(\mu_q\) is uniformly finite.
\end{proof}

\begin{lemma}\label{lem5.3}
For any $2< q< \frac{2n}{n-1}$ ($n\geq2$), there exists a quasi-positive function $\varphi_q$ satisfying \begin{equation}\label{5.5}
\square_{\omega_B}\varphi_q=N_1\mu_q\varphi^{q-1}_q
\end{equation} weakly in $W^{1,2}(M,\omega_B)$ with
\[
Y_q(\varphi_q)=\mu_q\quad\text{and}\quad\int_M\varphi_q^q\frac{\omega^n_B}{n!}=\frac{1}{N_1}.
\]
\end{lemma}

\begin{proof}
Fix $q\in(2,\frac{2n}{n-1})$, since \(Y_q(|\varphi|)= Y_q(\varphi)\) for every \(0\not\equiv\varphi\in W^{1,2}(M,\mathbb R)\), we may choose a minimizing real-valued
sequence \(\{\phi_i\}_{i\in\mathbb N}\subset W^{1,2}(M,\omega_B)\) such
that
\[
\phi_i\geq0,\qquad
Y_q(\phi_i)<\mu_q+\frac{1}{i},\qquad
\int_M\phi_i^q\frac{\omega_B^n}{n!}=\frac1{N_1}.
\]

By H\"older's inequality, we have
\[
\|\phi_i\|_2^2\leq \|\phi_i\|_q^2\leq N_1^{-\frac{2}{q}}.
\]
Moreover, by \eqref{5.3x} and \eqref{eq4.7}, we obtain
\begin{eqnarray*}
\|\partial\phi_i\|_2^2&\leq&Y_q(\phi_i)-N_1\int_MS_B^{(2)}(\omega_B)\phi_i^2\frac{\omega_B^n}{n!}\notag\\
&\leq&\big(1+N_1^{1-\frac{2}{q}}\big)\sup_M|S_B^{(2)}(\omega_B)|+1.
\end{eqnarray*}
Therefore, $\{\phi_i\}$ is bounded in $W^{1,2}(M,\omega_B)$.

We claim that there exists a quasi-positive function $\varphi_q\in W^{1,2}(M,\omega_B)$ satisfying
\[
Y_q(\varphi_q)=\mu_q\quad\text{and}\quad \int_M\varphi_q^q\frac{\omega_B^n}{n!}=\frac{1}{N_1}.
\]

Indeed, by the Rellich-Kondrakov embedding theorem (see, for instance, \cite[Theorem 2.34]{aubin}), $W^{1,2}(M,\omega_B)$ is compactly embedded in $L^q(M,\omega_B)$. Since $W^{1,2}(M,\omega_B)$ is reflexive, every bounded sequence admits a weakly convergent subsequence (see, for example, \cite[Appendix~D.4.]{evans}). Consequently, there exists a subsequence $\{\phi_{i_k}\}\subset\{\phi_i\}$ and $\varphi_q\in W^{1,2}(M,\omega_B)$ so that
\[
\text{(1)}\quad\phi_{i_k}\rightarrow \varphi_q\quad\text{strongly in}\quad L^q(M,\omega_B),
\]
and
\[
\text{(2)}\quad\phi_{i_k}\rightarrow \varphi_q\quad\text{weakly in}\quad W^{1,2}(M,\omega_B).
\]

It follows from (1) that
\begin{equation}\label{5.9x}
\int_M\varphi_q^q\frac{\omega_B^n}{n!}=\frac{1}{N_1}.
\end{equation}
and $\phi_{i_k}\rightarrow \varphi_q$ almost everywhere on $M$ (see, for instance, \cite[Proposition 3.43]{aubin}). Therefore, $\varphi_q$ is quasi-positive, i.e.,
\begin{equation}\label{5.10}
\varphi_q\geq0,\quad\text{and}\quad\varphi_q\not\equiv0.
\end{equation}

Since $q>2$, (1) and H\"older inequality shows that the convergence $\phi_{i_k}\rightarrow \varphi_q$ is also strong in $L^2(M,\omega_B)$. On the other hand, (2) implies that (see, for instance, \cite[Appendix~D.4]{evans})
\[\|\varphi_q\|_{1,2}\leq\liminf_{k\rightarrow\infty}\|\phi_{i_k}\|_{1,2}.\] Consequently, \[ \mu_q \leq Y_q(\varphi_q) \leq \liminf_{k\to\infty}Y_q(\phi_{i_k}) = \mu_q. \]
Hence $Y_q(\varphi_q)=\mu_q.$ This proves the claim.

Finally, we prove that $\varphi_q$ is a weak solution of \eqref{5.5} in  $W^{1,2}(M,\omega_B)$ using a variational method.

Since \(\varphi_q\) realizes the infimum of \(Y_q\) over
\(W^{1,2}(M,\mathbb R)\setminus\{0\}\), for every
\(\psi\in W^{1,2}(M,\omega_B)\) the first variation of
\(Y_q(\varphi_q+s\psi)\) at \(s=0\) vanishes.

Set $\varphi=\varphi_q+s\psi$ with $\psi\in W^{1,2}(M,\omega_B)$ and $s$ a small real parameter. It follows from \eqref{5.3x} and \eqref{5.9x} that
\begin{eqnarray}\label{5.11}
0&=&\frac{1}{2}\cdot\frac{d}{ds}\big|_{s=0}Y_q(\varphi)\notag\\
&=&\int_M\langle\partial\varphi_q,\partial\psi\rangle_{\omega_B}\frac{\omega_B^n}{n!}+N_1\int_MS_B^{(2)}(\omega_B)\varphi_q\psi\frac{\omega_B^n}{n!}\notag\\
&&-N_1\mu_q\int_M\varphi_q^{q-1}\psi\frac{\omega_B^n}{n!}.
\end{eqnarray}

Integrating by parts, and using \eqref{4.5} and $\eta(\omega_B)=0$, we obtain
\begin{equation}\label{5.12}
\int_M\langle\partial\varphi_q,\partial\psi\rangle_{\omega_B}\frac{\omega_B^n}{n!}=\frac{1}{2}\int_M\psi\Delta_d\varphi_q \frac{\omega_B^n}{n!}=-\int_M\psi\Delta_{\omega_B}^{\mathbb{C}}\varphi_q\frac{\omega_B^n}{n!}.
\end{equation}

Applying \eqref{5.12} to \eqref{5.11}, we get
\[\int_M\psi \square_{\omega_B}\varphi_q\frac{\omega_B^n}{n!}=N_1\mu_q\int_M\varphi_q^{q-1}\psi\frac{\omega_B^n}{n!}\]
for all $\psi\in W^{1,2}(M,\omega_B)$, i.e., $\varphi_q$ solves \eqref{5.5} weakly in $W^{1,2}(M,\omega_B)$.
\end{proof}

\begin{theorem}\label{thm5.1}
For $2< q< \frac{2n}{n-1}$ ($n\geq2$), there exists a strictly positive function $\varphi_q\in C^\infty(M,\mathbb{R})$ satisfying \eqref{5.5} with
\[Y_q(\varphi_q)=\mu_q.\]
\end{theorem}

\begin{proof}
Let \(\varphi_q\) be the nonnegative weak solution obtained in
Lemma \ref{lem5.3}, i.e.,
\begin{equation}\label{eq4.15}
\Delta_{\omega_B}^\mathbb{C}\varphi_q=-N_1\mu_q\varphi_q^{q-1}+N_1S_B^{(2)}(\omega_B)\varphi_q.
\end{equation}
Since \(\varphi_q\in W^{1,2}(M,\omega_B)\), the Sobolev embedding yields
\(
\varphi_q\in L^{\frac{2n}{n-1}}(M,\omega_B).
\)

We claim that \(\varphi_q\in L^\infty(M)\).

Set
\[
r_0:=\frac{2n}{n-1}.
\]
Suppose that \(\varphi_q\in L^{r_k}(M)\). Since \(q\in(2,r_0)\), the right-hand side of
\eqref{eq4.15} belongs to
\[
L^{p_k}(M),
\qquad
p_k:=\frac{r_k}{q-1}>1.
\]
Hence \(W^{2,p_k}\)-estimate (Calder\'on-Zygmund inequality) for elliptic equations implies that $\varphi_q\in W^{2,p_k}(M)$. If \(p_k<n\), then the Sobolev embedding in real dimension \(2n\) yields
\begin{equation}\label{eq4.16}
\varphi_q\in L^{r_{k+1}}(M),
\qquad
\frac{1}{r_{k+1}}=\frac{q-1}{r_k}-\frac{1}{n}.
\end{equation}

Writing \(x_k:=r_k^{-1}\), we have
\[
x_{k+1}=(q-1)x_k-\frac1n,\qquad
x_0=\frac{n-1}{2n}<\frac{1}{n(q-2)},
\]
The latter inequality is equivalent to
\(
q<\frac{2n}{n-1}.
\)
Moreover,
\(
p_0=\frac{r_0}{q-1}>1.
\)
As long as \(p_k<n\), the recursion above gives
\[
x_{k+1}-x_k=(q-2)x_k-\frac1n<0,
\]
and hence \(r_{k+1}>r_k\). In particular, \(p_k>1\) throughout the
iteration.

Assume that \(p_k<n\) for every \(k\). Then \eqref{eq4.16} is valid for
all \(k\), and
\[
x_k-\frac{1}{n(q-2)}
=(q-1)^k
\left(
x_0-\frac{1}{n(q-2)}
\right).
\]
Hence \(x_k<0\) for sufficiently large \(k\), which is impossible.
Therefore, \(p_k\geq n\) for some \(k\).

If \(p_k>n\), then the Sobolev embedding
\(
W^{2,p_k}(M)\hookrightarrow L^\infty(M)
\)
implies that \(\varphi_q\in L^\infty(M)\). If \(p_k=n\), then
\(
W^{2,n}(M)\hookrightarrow L^r(M)
\)
for every \(r<\infty\). Choosing \(r>n(q-1)\), we obtain
\(\varphi_q\in L^r(M)\), and therefore
\[
\varphi_q^{q-1}\in L^{\frac{r}{q-1}}(M),
\qquad
\frac{r}{q-1}>n.
\]
Since \(S_B^{(2)}(\omega_B)\varphi_q\in L^r(M)\), it also belongs to
\(L^{\frac{r}{q-1}}(M)\). Applying the \(W^{2,p}\)-estimate to
\eqref{eq4.15}, with
\(
p:=\frac{r}{q-1}>n,
\)
we obtain \(\varphi_q\in W^{2,p}(M)\), and hence
\[
\varphi_q\in L^\infty(M).
\]

It follows that the right-hand side of \eqref{eq4.15},
\[
a(x):=
-N_1S_B^{(2)}(\omega_B)
+
N_1\mu_q\varphi_q^{q-2},
\]
belongs to
\(L^\infty(M)\). Therefore,  \(W^{2,p}\)-estimates give
\(
\varphi_q\in W^{2,p}(M)
\)
for every \(p<\infty\). Taking \(p>2n\),  the Rellich-Kondrakov embedding theorem (see, for instance, \cite[Theorem 2.34]{aubin}) implies
\(
\varphi_q\in C^{1,\gamma}(M)
\)
for some \(\gamma\in(0,1)\). 

Choose
\[
\Lambda:=\|a\|_{\infty}+1.
\]
Then \eqref{eq4.15} gives
\[
\bigl(\Delta_{\omega_B}^C-\Lambda\bigr)(-\varphi_q)\geq0.
\]
If \(\varphi_q\) vanishes at some point of \(M\), then
\(-\varphi_q\) attains its nonnegative maximum \(0\) there and  the strong maximum principle
implies that \(\varphi_q\equiv0\), contradicting
\[
\int_M\varphi_q^q\frac{\omega_B^n}{n!}=\frac1{N_1}.
\]
Therefore,
\(
\varphi_q>0
\)
on \(M\).

Since \(M\) is compact, there exists
\(c>0\) such that \(\varphi_q\geq c\). The map
\(s\mapsto s^{q-1}\) is smooth on \([c,\infty)\). Hence \eqref{eq4.15},
together with Schauder estimates and elliptic bootstrapping, yields
\[
\varphi_q\in C^\infty(M,\mathbb R).
\]
\end{proof}

Now we are ready to complete the proof of the first main theorem.\par
\medskip

\noindent\textbf{Proof of Theorem \ref{thm1}.} Since $n\geq2$, we have $2<N_2<\frac{2n}{n-1}$. By Theorem \ref{thm5.1}, there exists a smooth and strictly positive solution $\varphi_{N_2}$ to \eqref{5.5}  with $q=N_2$, satisfying
$$Y_{N_2}(\varphi_{N_2})=\mu_{N_2}.$$

Set $\lambda:=\mu_{N_2}$, then \eqref{5.5} is precisely \eqref{5.2}. It follows that $$f:=\frac{2n-1}{n^2-1}\log\varphi_{N_2}\in C^\infty(M,\mathbb{R})$$
solving \eqref{5.1}. Therefore
$$\widetilde{\omega}=e^{f}\omega_G\in\{\omega\}$$ has constant second Bismut scalar curvature $$S_B^{(2)}(\widetilde{\omega})=\inf_{\varphi\in C^\infty(M,\mathbb{R})}Y_{N_2}(\varphi).$$$\hfill\Box$

\medskip

\noindent\textbf{Proof of Corollary \ref{cor4}.} Since \(\omega_G\) is Gauduchon, we have
\(d^*\eta(\omega_G)=0\). Together with $d\eta(\omega_G)=0$, \(\eta(\omega_G)\) is a harmonic $1$-form.
By the Hodge theorem,
\[
\dim\mathcal H^1(M)=b_1(M)=0,
\]
and therefore \(\eta(\omega_G)=0\). Thus
\(\omega_G\) is balanced, so the conformal class \(\{\omega\}\) is
balanced. The conclusion follows from Theorem \ref{thm1}.
\section{Second Chern scalar curvature and geometric consequences}
\label{sec:5}

In this section, we study the existence of constant second Chern scalar curvature within a Hermitian conformal class on compact Hermitian manifolds and derive its geometric consequences.

The existence proofs follow the strategy of \cite[Theorems 3.1, 4.1]{ACS17} for the Chern--Yamabe problem (see also \cite{Bar23} for the Bismut analogue), with minor modifications to account for the Laplacian coefficient difference.

We will provide detailed proofs of Theorems~\ref{thm2} and \ref{thm3}, covering both the existence results and the geometric applications.\par
\medskip

\noindent\textbf{Proof of Theorem \ref{thm2}.} By \eqref{4.2}, the existence problem reduces to solving
\begin{equation}\label{6.1}
\Delta^{\mathbb{C}}_{\omega_G} f=S_C^{(2)}(\omega_G).
\end{equation}

Combining \cite[Lemma 2.2]{YZ25} and \eqref{4.5}, the formal adjoint $(\Delta_{\omega_G}^\mathbb{C})^*$ of $\Delta_{\omega_G}^\mathbb{C}$ satisfies
\begin{equation}\label{6.2}
(\Delta_{\omega_G}^\mathbb{C})^*=-\frac{1}{2}\Delta_d+\frac{1}{2}\langle d(\cdot),\eta(\omega_G)\rangle_{\omega_G}.
\end{equation}

We claim that $\mathrm{Ker}(\Delta_{\omega_G}^\mathbb{C})^*$ consists of constant functions.

Indeed, for any $u \in \mathrm{Ker}(\Delta_{\omega_G}^\mathbb{C})^*$, since $\omega_G$ is Gauduchon, we have
\begin{eqnarray*}
0&=&(u,(\Delta_{\omega_G}^\mathbb{C})^*u)_{\omega_G}\notag\\
&=&-\frac{1}{2}(u,\Delta_du)_{\omega_G}+\frac{1}{4}(u^2,d^*\eta(\omega_G))_{\omega_G}\notag\\
&=&-\frac{1}{2}\|du\|^2_{\omega_G}.
\end{eqnarray*}
Hence $du=0$, and $u$ is constant.

Together with the assumption
\[\Gamma_M^{(2)}(\{\omega\})=\int_MS_C^{(2)}(\omega_G)\frac{\omega_G^n}{n!}=0,\]
we obtain
\[S_C^{(2)}(\omega_G)\in (\mathrm{Ker}(\Delta_{\omega_G}^\mathbb{C})^*)^\perp=\mathrm{Im}\Delta_{\omega_G}^\mathbb{C}.\]

Therefore, \eqref{6.1} has a unique solution $f\in C^\infty(M,\mathbb{R})$, up to additive constants. This proves the first assertion of the theorem, namely, that there exists a unique (up to scaling) Hermitian metric
\[\omega_a:=e^f\omega_G\in\{\omega\}\]
with $S_C^{(2)}(\omega_a)=0$.

Moreover, by taking $t=0$ in \eqref{7.1} gives that
\begin{equation}\label{6.3}
\Gamma_M^{(1)}(\{\omega\})=\int_MS^{(1)}_C(\omega_G)\frac{\omega^n_G}{n!}=\|\bar{\partial}^*\omega_G\|^2_{\omega_G}\geq0.
\end{equation}

If $\Gamma_M^{(1)}(\{\omega\})>0$, then $\kappa(M)=-\infty$ by \cite[Corollary 3.3]{Yang19}, while if $\Gamma_M^{(1)}(\{\omega\})=0$, then \cite[Theorem 1.4]{Yang19} shows that $M$ lies in one of the following cases:

\begin{enumerate}
\item $\kappa(M)=-\infty$ and neither $K_M$ nor $K_M^{-1}$ is pseudo-effective;
\item $\kappa(M)=-\infty$ and $K_M$ is unitary flat;
\item $\kappa(M)=0$ and $K_M$ is a holomorphic torsion, i.e., $K_M^{\otimes m}=\mathcal{O}_M$ for some $m\in\mathbb{Z}^+$.
\end{enumerate}
This proves the second assertion of the theorem, namely, that either $\kappa(M)=-\infty$, or $\kappa(M)=0$ and $K_M$ is a holomorphic torsion.

If, in addition, $\{\omega\}$ is balanced, then \eqref{7.1} implies that $S_C^{(1)}(\omega_G)=S_C^{(2)}(\omega_G)$. By \cite[Theorem 3.1]{ACS17}, there exists a unique (up to scaling) metric $\omega_b\in\{\omega\}$ such that $S_C^{(1)}(\omega_b)=0$, proving the last assertion of the theorem.$\hfill\Box$\par
\medskip

We divide the proof of Theorem \ref{thm3} into several lemmas.

\begin{lemma}\label{lem6.1}
Let $(M,\omega)$ be a compact Hermitian manifold with complex dimension $n\geq2$. If $\Gamma_M^{(2)}(\{\omega\})<0$, then there exists a representative $\omega' \in \{\omega\}$ whose second Chern scalar curvature is pointwise negative.
\end{lemma}

\begin{proof}
It is well known (see, for instance, \cite{Gau77,Gau84}) that there exists a unique (up to additive constants) solution $u\in C^{\infty}(M,\mathbb{R})$ of
\begin{equation}\label{6.4}
\Delta_{\omega_G}^{\mathbb{C}}u=S_C^{(2)}(\omega_G)-\frac{\int_MS_C^{(2)}(\omega_G)\omega_G^n}{\int_M\omega_G^n}
\end{equation}
since $\omega_G$ is Gauduchon and
\[\int_M\big(S_C^{(2)}(\omega_G)-\frac{\int_MS_C^{(2)}(\omega_G)\omega_G^n}{\int_M\omega_G^n}\big)\frac{\omega_G^n}{n!}=0.\]

Consider $\omega'=e^u\omega_G\in\{\omega\}$. It follows from \eqref{4.2}, \eqref{6.4} and $\mathrm{Vol}(M,\omega_G)=1$ that
\begin{eqnarray}\label{6.5}
S_C^{(2)}(\omega')&=&e^{-u}(S_C^{(2)}(\omega_G)-\Delta_{\omega_G}^{\mathbb{C}}u)\notag\\
&=&e^{-u}\cdot\frac{\int_MS_C^{(2)}(\omega_G)\omega_G^n}{\int_M\omega_G^n}\notag\\
&=&e^{-u}\Gamma_M^{(2)}(\{\omega\})<0.
\end{eqnarray}
\end{proof}

The Yamabe-type problem of finding a function $f\in C^\infty(M,\mathbb{R})$ such that the conformally related metric
\[
\omega_b := e^{f}\omega'
\]
has negative constant second Chern scalar curvature that
\[
S_C^{(2)}(\omega_b)=\Gamma_M^{(2)}(\{\omega\})<0
\]
is equivalent to proving the existence of a smooth solution to the equation
\begin{equation}\label{6.6}
\Delta_{\omega'}^{\mathbb{C}} f = -\lambda e^{f} + S_C^{(2)}(\omega')\quad\text{with}\quad \lambda=\Gamma_M^{(2)}(\{\omega\}).
\end{equation}
We shall solve \eqref{6.6} by the continuity method.

For $\alpha\in(0,1)$, we define the map
\[F:\ [0,1]\times C^{2,\alpha}(M,\mathbb{R})\rightarrow C^{\alpha}(M,\mathbb{R})\]
by
\begin{equation}\label{6.7}
F(a,f)=\Delta_{\omega'}^{\mathbb{C}}f-aS_C^{(2)}(\omega')+\lambda e^f-\lambda(1-a),
\end{equation}
and define the set
\begin{equation}\label{6.8}
A := \Big\{ a\in[0,1]\;\Big|\; F(a,f_a)=0 \ \text{for some } f_a\in C^{2,\alpha}(M,\mathbb{R}) \Big\}.
\end{equation}

\begin{lemma}\label{lem6.2}
The set $A$ defined in \eqref{6.8} is open in $[0,1]$.
\end{lemma}

\begin{proof}
Since $F(0,0)=0$, $A$ is nonempty.

Fix $a_0\in A$, with corresponding $f_{a_0}\in C^{2,\alpha}(M,\mathbb{R})$ such that $F(a_0,f_{a_0})=0$. The linearization of $F$ with respect to $f$ at
$(a_0,f_{a_0})$ is the operator
\[D:\, C^{2,\alpha}(M,\mathbb{R})\rightarrow C^{\alpha}(M,\mathbb{R}),\qquad Dw:=\Delta_{\omega'}^{\mathbb{C}}w+\lambda e^{f_{a_0}}w.\]
Clearly, $D$ is a second-order linear elliptic operator.

Suppose that $w\in\mathrm{Ker}D$. Let $p,\ q\in M$ be points at which $w$ attains its maximum and minimum, respectively. Combining with the maximum principle, we have
\[\lambda e^{f_{a_0}(p)}w(p)\geq0\quad\text{and}\quad \lambda e^{f_{a_0}(q)}w(q)\leq0.\]
Since $\lambda<0$, it follows that
\[0\geq w(p)\geq w(q)\geq0,\]
hence $w=0$ on $M$. Therefore, $\mathrm{Ker}D = \{0\}$ and $D$ is injective.

Since \(\Delta_{\omega'}^{\mathbb C}\) has index zero and \(D\) differs from it by a zeroth-order term, \(D\) is a Fredholm operator of index zero.
It follows that
\[
\dim\operatorname{Coker}D=\dim\operatorname{Ker}D-\operatorname{ind}D=0.
\]
Hence \(D\) is surjective.

By the implicit function theorem in Banach spaces, the set $A \subset [0,1]$ is open.
\end{proof}

We now turn to proving that $A$ is also closed.

For any $\{a_n\}\subset A$ and corresponding functions $f_{a_n}\in C^{2,\alpha}(M,\mathbb{R})$ satisfy $F(a_n,f_{a_n})=0$ for each $n$. It follows from \eqref{6.7} that $f_{a_n}$ satisfies
\begin{equation}\label{6.12}
\Delta_{\omega'}^{\mathbb{C}} f_{a_n} = a_n S_C^{(2)}(\omega') -\lambda e^{f_{a_n}} + \lambda(1-a_n).
\end{equation}

To apply the Arzel\'{a}--Ascoli theorem and extract a subsequence converging in $C^{2,\alpha}(M,\mathbb{R})$ to a limit $f_{a_\infty}$ with $F(a_\infty,f_{a_\infty})=0$, we need to establish a series of a priori estimates for the functions $f_{a_n}$.

\begin{lemma}\label{lem6.3}
For all $n$ and all $x\in M$, the functions $f_{a_n}$ satisfy
\begin{equation}\label{6.9}
\min\big\{0,\log\frac{\max_MS_C^{(2)}(\omega')}{\lambda}\big\}\leq f_{a_n}(x)\leq \log\big(1+\frac{\min_MS_C^{(2)}(\omega')}{\lambda}\big).
\end{equation}
\end{lemma}

\begin{proof}
Fix $n$, and let $p\in M$ (resp. $q\in M$) be a maximum (resp. minimum) point of $f_{a_n}$.

Using \eqref{6.12}, we obtain that at $p$,
\[
a_n S_C^{(2)}(\omega')-\lambda e^{f_{a_n}}+\lambda(1-a_n)\leq0.
\]
Since $\lambda<0$, $S_C^{(2)}(\omega')<0$ and $0\leq a_n\leq 1$, taking logarithms yields
\begin{equation}\label{6.10}
f_{a_n}(x)\leq f_{a_n}(p)\leq\log\big(1+\frac{\min_MS_C^{(2)}(\omega')}{\lambda}\big)
\end{equation}
for all $x\in M$.

Similarly, at $q$ we have
\[
a_n S_C^{(2)}(\omega')-\lambda e^{f_{a_n}}+\lambda(1-a_n)\geq0,
\]
which implies
\begin{equation}\label{6.11}
f_{a_n}(x)\geq f_{a_n}(q)\geq\min\big\{0,\log\frac{\max_MS_C^{(2)}(\omega')}{\lambda}\big\}.
\end{equation}
for all $x\in M$.

Combining \eqref{6.10} and \eqref{6.11} gives \eqref{6.9}.
\end{proof}

\begin{lemma}\label{lem6.4}
The set $A$ defined in \eqref{6.8} is closed in $[0,1]$.
\end{lemma}

\begin{proof}
Let $\{a_n\}\subset A$ with $a_n\to a_\infty\in[0,1]$, and let
$f_{a_n}\in C^{2,\alpha}(M)$ satisfy \eqref{6.12}.
By Lemma~\ref{lem6.3}, the sequence $\{f_{a_n}\}$ is uniformly bounded in $C^0(M)$. Hence, the right-hand-side of \eqref{6.12} is uniformly bounded in $L^\infty(M,\omega')$, and therefore in $L^p(M,\omega')$ for every $p>1$.

Since $\Delta_{\omega'}^{\mathbb{C}}$ is uniformly elliptic,
the Calder\'on-Zygmund inequality yields uniform $W^{2,p}(M,\omega')$ bounds for any $p>1$.
Choose $p>\frac{2n}{1-\alpha}$, and set $\gamma:=1-\frac{2n}{p}>\alpha$.
By Sobolev embedding in real dimension \(2n\), the sequence
\(\{f_{a_n}\}\) is uniformly bounded in \(C^{1,\gamma}(M)\). It follows from \eqref{6.12} that
\({\Delta_{\omega_0}^{\mathbb C}f_{a_n}}\) is uniformly bounded in
\(C^\alpha(M,\mathbb{R})\). Schauder estimates then imply that
\(\{f_{a_n}\}\) is uniformly bounded in \(C^{2,\alpha}(M,\mathbb{R})\).

By the Ascoli-Arzel\'{a} theorem, after passing to a subsequence,
$f_{a_n}\to f_{a_\infty}$ in $C^{2,\beta}(M)$ for any $\beta<\alpha$.
Passing to the limit in \eqref{6.12} shows
that \(f_{a_\infty}\) is a \(C^{2,\beta}\) solution of \eqref{6.12} with
\(a=a_\infty\). Standard elliptic bootstrapping yields
\(f_{a_\infty}\in C^\infty(M,\mathbb{R})\). Hence \(a_\infty\in A\), and therefore
\(A\) is closed.
\end{proof}

\noindent\textbf{Proof of Theorem \ref{thm3}.} Combining the openness and closedness of the set $A$, the continuity method yields the existence of a $C^{2,\alpha}$ solution to \eqref{6.6}. By standard elliptic regularity, this solution is smooth.

We next prove uniqueness. Suppose that both $f_1$ and $f_2\in C^\infty(M,\mathbb{R})$ solving \eqref{6.6}, namely,
\[
\Delta_{\omega'}^{\mathbb{C}} f_1 = -\lambda e^{f_1} + S_C^{(2)}(\omega'),\qquad \Delta_{\omega'}^{\mathbb{C}} f_2 = -\lambda e^{f_2} + S_C^{(2)}(\omega').
\]
Subtracting the two equations gives
\begin{equation}\label{6.13}
\Delta_{\omega'}^{\mathbb{C}} (f_1-f_2) = \lambda ( e^{f_2}-e^{f_1}).
\end{equation}

Let $p,\ q\in M$ be points where $f_1-f_2$ attains its maximum and minimum, respectively. Evaluating \eqref{6.13} at these points and applying the maximum principle, we obtain
\[
0\leq (f_1-f_2)(q)\leq (f_1-f_2)(p)\leq 0,
\]
which forces $f_1=f_2$ on $M$. Hence the solution of \eqref{6.6} is unique.

This proves the first assertion of the theorem: there exists a unique (up to scaling) Hermitian metric
\[\omega_c:=e^{f+u}\omega_G\in\{\omega\}\]
satisfying
\[
S_C^{(2)}(\omega_c)=\Gamma_M^{(2)}(\{\omega\})<0.
\]

If, in addition, $\{\omega\}$ is balanced, i.e., $\omega_G$ is balanced, then by taking $t=0$ in \eqref{7.1} we have $S_C^{(1)}(\omega_G)=S_C^{(2)}(\omega_G)$ . Consequently,
\[\Gamma_M^{(1)}(\{\omega\})=\Gamma_M^{(2)}(\{\omega\})<0.\]

By \cite[Theorem 1.3]{Yang19}, this implies that $K_M^{-1}$ is not pseudo-effective, proving the second assertion of the theorem.

Finally, \cite[Theorem 4.1]{ACS17} guarantees the existence of a unique (up to scaling) Hermitian metric \(\omega_d \in \{\omega\}\) such that
\[
S_C^{(1)}(\omega_d) = \Gamma_M^{(2)}(\{\omega\}) < 0.
\]
$\hfill\Box$\par
\medskip

We now turn to the complementary case $\Gamma_M^{(2)}(\{\omega\})>0$, which leads to a substantially different geometric consequence.\par
\medskip

\begin{proposition}\label{thm4}
Let $(M,\omega)$ be a compact Hermitian manifold of complex dimension $n\ge 2$ and let $\mathcal{C}$ be a compact complex curve of genus $g(\mathcal{C})\ge 2$.
If $\Gamma_M^{(2)}(\{\omega\}) > 0$, then $\kappa(M)=-\infty$ and there exists a Hermitian metric on $M\times \mathcal{C}$ with positive constant first Chern scalar curvature, as well as Hermitian metrics with positive constant second Chern scalar curvature.
\end{proposition}

\begin{proof}
If $\Gamma_M^{(2)}(\{\omega\})>0$, then taking $t=0$ in \eqref{7.1} yields $\Gamma_M^{(1)}(\{\omega\})>0$.
By \cite[Corollary 3.3]{Yang19}, this implies that $\kappa(M)=-\infty$.

The existence of a Hermitian metric with positive constant first Chern scalar curvature follows from \cite[Proposition 5.7]{ACS17}.
In addition, the proof of \cite[Lemma 5.8]{ACS17} applies verbatim in the case $\Gamma_M^{(2)}(\{\omega\})>0$, yielding the corresponding result for the second Chern scalar curvature.
\end{proof}

\section{Weak second Hermitian--Einstein metrics}
\label{sec:6}

In this section, we study weak second Hermitian--Einstein metrics. We prove Theorem~\ref{thm5}, showing that on compact pluriclosed Gauduchon manifolds such a condition implies the constancy of the second Chern scalar curvature, and derive the corresponding corollaries \ref{cor1}, \ref{cor2} and \ref{cor3}.\par
\medskip

\noindent\textbf{Proof of Theorem \ref{thm5}.} By setting $t=0$ in \eqref{1.2b} and \eqref{1.2c}, we obtain
\begin{equation}\label{5.20x}
\Theta^{(3)}(\omega)+\Theta^{(4)}(\omega)=2\Theta^{(1)}-(\partial\partial^*\omega+\bar{\partial}\bar{\partial}^*\omega).
\end{equation}
In particular,
\begin{equation}\label{5.20}
\partial\bar{\partial}(\Theta^{(3)}(\omega)+\Theta^{(4)}(\omega))=0.
\end{equation}

On the other side, since $\omega$ is pluriclosed, namely $\partial\bar{\partial}\omega=0$, a direct computation yields
\begin{eqnarray}\label{5.21}
&&\partial\bar{\partial}(f\omega)\notag\\
&=&\sqrt{-1}\frac{\partial^2}{\partial z^i\partial\bar{z}^j}(fh_{k\bar{l}})dz^i\wedge d\bar{z}^j\wedge dz^k\wedge d\bar{z}^l\notag\\
&=&\sqrt{-1}\big(\frac{\partial^2f}{\partial z^i\partial\bar{z}^j}h_{k\bar{l}}+\frac{\partial f}{\partial z^i}\frac{\partial h_{k\bar{l}}}{\partial\bar{z}^j}+\frac{\partial f}{\partial \bar{z}^j}\frac{\partial h_{k\bar{l}}}{\partial z^i}\big)dz^i\wedge d\bar{z}^j\wedge dz^k\wedge d\bar{z}^l.
\end{eqnarray}

Contracting \eqref{5.21} with
\(
h^{i\bar j}h^{k\bar{l}}-h^{i\bar{l}}h^{k\bar j},
\)
we obtain
\begin{eqnarray}\label{5.22}
&&\left(h^{i\bar j}h^{k\bar\ell}
-h^{i\bar\ell}h^{k\bar j}\right)
\left(
\frac{\partial^2 f}{\partial z^i\partial\bar z^j}h_{k\bar\ell}
+\frac{\partial f}{\partial z^i}
\frac{\partial h_{k\bar\ell}}{\partial\bar z^j}
+\frac{\partial f}{\partial\bar z^j}
\frac{\partial h_{k\bar\ell}}{\partial z^i}
\right) \notag\\
&=&
(n-1)\Delta_\omega^{\mathbb C}f
+
2\operatorname{Re}
\left\langle
\sqrt{-1}\partial f,\bar\partial^\ast\omega
\right\rangle_\omega .
\end{eqnarray}

Combining \eqref{2CE}, \eqref{5.20}, \eqref{5.21}, and \eqref{5.22}, we conclude that $f$ satisfies the linear elliptic equation
\begin{equation}\label{5.23x}
\lozenge_\omega f:=(n-1)\Delta_\omega^\mathbb{C}f+2\mathrm{Re}\langle\sqrt{-1}\partial f,\bar{\partial}^*\omega\rangle_\omega=0.
\end{equation}

Since $\omega$ is Gauduchon,
\[
d^*\eta(\omega)=0,\quad\text{and}\quad \bar{\partial}^*\partial^*\omega=0.
\]
For any $u\in \mathrm{Ker}\lozenge_\omega$, integrating by parts and using \eqref{4.5} gives that
\begin{eqnarray*}
0&=&\int_M u\lozenge_\omega u\frac{\omega^n}{n!}\notag\\
&=&(n-1)\int_Mu\Delta_\omega^\mathbb{C}u\frac{\omega^n}{n!}+2\int_Mu\mathrm{Re}\langle\sqrt{-1}\partial u,\bar{\partial}^*\omega\rangle_\omega\frac{\omega^n}{n!}\notag\\
&=&-\frac{n-1}{2}\int_Mu\Delta_du-\frac{n-1}{4}\int_M\langle du^2,\eta(\omega)\rangle_\omega+\int_M\mathrm{Re}\langle\sqrt{-1}\partial u^2,\bar{\partial}^*\omega\rangle_\omega\frac{\omega^n}{n!}\notag\\
&=&-\frac{n-1}{2}\|du\|^2_\omega-\frac{n-1}{4}(u^2,d^*\eta(\omega))_\omega+(\sqrt{-1}\bar{\partial}^*\partial^*\omega, u^2)_\omega\notag\\
&=&-\frac{n-1}{2}\|du\|^2_\omega,
\end{eqnarray*}
It follows that $u$ is a constant, and therefore all solutions of \eqref{5.23x} are constants.

Finally, tracing \eqref{2CE}, we conclude that
\[S_C^{(2)}=\frac{n}{2}f=\mathrm{const}.\]

This completes the proof.$\hfill\Box$\par
\medskip

We now apply Theorem~\ref{thm5} to prove Corollaries \ref{cor1} and \ref{cor2}.\par
\medskip

\noindent\textbf{Proof of Corollary \ref{cor1}.} It follows from \eqref{5.20x} and \cite[(4.9)]{Yang25b+} that
\begin{eqnarray}\label{5.23}
(\partial(\Theta^{(3)}(\omega)+\Theta^{(4)}(\omega)),\partial\omega)_\omega&=&(-\partial\bar{\partial}\bar{\partial}^*\omega,\partial\omega)_\omega\notag\\
&=&-(\bar{\partial}\bar{\partial}^*\omega,\partial^*\partial\omega)_\omega\notag\\
&=&\|\partial\bar{\partial}^*\omega\|^2_\omega.
\end{eqnarray}

In complex dimension two, any Gauduchon metric is also pluriclosed. Hence, by Theorem~\ref{thm5}, the function $f=\frac{2}{n}S_C^{(2)}(\omega)$ is constant.

Combining \eqref{2CE} and \eqref{5.23}, we have
\begin{equation}\label{5.24}
\|\partial\bar{\partial}^*\omega\|^2_\omega=f\|\partial\omega\|_\omega^2.
\end{equation}

If $f=0$, then $S_C^{(2)}(\omega)=0$.

If $f< 0$, \eqref{5.24} forces $\partial\omega=0$, hence $\omega$ is K\"ahler. In this case the first to fourth Chern--Ricci curvatures coincide, and \eqref{2CE} gives
\[Ric(\omega)=\frac{f}{2}\omega.\]
Thus $\omega$ is K\"ahler--Einstein with negative scalar curvature.
$\hfill\Box$\par
\medskip

\noindent\textbf{Proof of Corollary \ref{cor2}.} Since $\omega$ is balanced, it is in particular Gauduchon and satisfies
\begin{equation}\label{5.25}
d(\Theta^{(3)}(\omega)+\Theta^{(4)}(\omega))=2d\Theta^{(1)}(\omega)=0.
\end{equation}

By Theorem~\ref{thm5}, the function $f$ is constant. There are two cases:

If $f=0$, then $S_C^{(2)}(\omega)\equiv0$.

If $f\neq0$, then \eqref{2CE} and \eqref{5.25} imply $d\omega = 0$, so that $\omega$ is K\"ahler--Einstein with non-zero scalar curvature.$\hfill\Box$\par
\medskip

We now prove Corollary \ref{cor3} as a consequence of Theorems \ref{thm2} and \ref{thm3}.\par
\medskip

\noindent\textbf{Proof of Corollary \ref{cor3}.} Let $\omega_G:=e^u\omega$ be the volume-normalized Gauduchon representative of $\{\omega\}$.

By \eqref{2CE} and \eqref{4.2}, we obtain
\begin{eqnarray}\label{5.26}
S_C^{(2)}(\omega_G)&=&e^{-u}\big(S_C^{(2)}(\omega)-\Delta_{\omega}^{\mathbb{C}}u\big)\notag\\
&=&\frac{n}{2}e^{-u}f-\Delta_{\omega_G}^{\mathbb{C}}u.
\end{eqnarray}

Note that $f\leq0$ by assumption. Together with \eqref{4.5}, \eqref{5.26} implies that
\begin{eqnarray}
\Gamma_M^{(2)}(\{\omega\})&=&\int_M S_C^{(2)}(\omega_G)\frac{\omega_G^n}{n!}\notag\\
&=&\frac{n}{2}(e^{-u},f)_{\omega_G}+\frac{1}{2}(\Delta_du,1)_{\omega_G}+\frac{1}{2}(u,d^*\eta(\omega_G))_{\omega_G}\notag\\
&=&\frac{n}{2}(e^{-u},f)_{\omega_G}\leq0.
\end{eqnarray}
If $\Gamma_M^{(2)}(\{\omega\})=0$, the conclusion follows from Theorem \ref{thm2}, while if $\Gamma_M^{(2)}(\{\omega\})<0$, it follows from Theorem \ref{thm3}.$\hfill\Box$

\section{Non-K\"ahler examples of constant Chern scalar curvatures}
\label{sec:7}
To illustrate the main results of the previous sections, we present explicit examples of Hermitian metrics with constant second Chern scalar curvature.

\noindent\textbf{Example 7.1.} Let $\mathbb{S}^{2n-1}\times\mathbb{S}^1$ be the standard $n$-dimensional ($n\geq2$) Hopf manifold. It is diffeomorphic to $(\mathbb{C}^n)^*/G$, where $(\mathbb{C}^n)^*=\mathbb{C}^n\backslash\{0\}$ and $G$ is a cyclic group generated by the contraction $z\mapsto \frac{1}{2}z$. This manifold carries a natural complex structure and a Hermitian metric
\[
\omega_h:=\sqrt{-1}\, \frac{4\delta_{ij}}{|z|^2}dz^i\wedge d\bar{z}^j
\]
which is induced from $(\mathbb{C}^n\backslash\{0\},\omega_{\mathrm{can}})$.

It was computed in \cite{LY12,LY17} that
\[\Theta(\omega_h)_{i\bar{j}k\bar{l}}=\frac{4\delta_{kl}(\delta_{ij}|z|^2-z^i\bar{z}^j)}{|z|^6},\]
\[\Theta^{(1)}(\omega_h)=\sqrt{-1}\, \partial\bar{\partial}\log|z|^{2n}\geq0,\qquad \Theta^{(2)}(\omega_h)=\frac{n-1}{4}\omega_h>0,\]
and
\[S_C^{(1)}(\omega_h)=\frac{n(n-1)}{4}.\]

It is straightforward to verify that
\begin{equation}\label{8.1}
\Theta^{(3)}(\omega_h)=\Theta^{(4)}(\omega_h)=\sqrt{-1}\, \frac{\delta_{ij}|z|^2-z^i\bar{z}^j}{|z|^4}dz^i\wedge d\bar{z}^j\geq0,
\end{equation}
and
\begin{equation}\label{8.2}
S_C^{(2)}(\omega_h)=\frac{n-1}{4}.
\end{equation}

\noindent\textbf{Example 7.2.} We consider a non-K\"ahler properly elliptic surface drawn from Tosatti and Weinkove \cite[Section 8]{TW13} (see also Yang \cite{Yang25}). It is constructed as follows.

Let $\pi:\, S\rightarrow C$ be an elliptic fiber bundle, where $S$ is a minimal non-K\"ahler compact complex surface of Kodaira dimension $\kappa(S)=1$ and the base curve $C$ has genus $g(C)\geq2$. The universal cover of $S$ can be identified biholomorphically with  $\mathbb{H}\times\mathbb{C}$ where $\mathbb{H}:=\{z\in\mathbb{C}\, |\, \mathrm{Im}z>0\}$ denotes the upper half plane. Define a holomorphic covering map
\[
h:\, \mathbb{H}\times\mathbb{C}\rightarrow \mathbb{H}\times\mathbb{C}^*,\qquad h(z,z')=(z,w)=(z,e^{-\frac{z'}{2}}),
\]
where $z=x+\sqrt{-1}\, y$.

Vaisman \cite{Vai87} discovered a Gauduchon metric on $\mathbb{H}\times\mathbb{C}^*$ given by
\[
\omega_0=\sqrt{-1}\, \big(-\frac{2}{w}dw+\frac{\sqrt{-1}}{y}dz\big)\wedge\big(-\frac{2}{\bar{w}}d\bar{w}-\frac{\sqrt{-1}}{y}d\bar{z}\big)+\frac{1}{y^2}\sqrt{-1}\, dz\wedge d\bar{z}
\]
which satisfies
\[
\bar{\partial}\omega_0=-\frac{\sqrt{-1}}{y^2\bar{w}}d\bar{z}\wedge dz\wedge d\bar{w}.
\]
$\omega_0$ descends to a smooth Hermitian metric on $S$.

Tosatti and Weinkove \cite{TW13} computed that
\begin{equation}\label{8.3}
\Theta^{(1)}(\omega_0)=-\frac{\sqrt{-1}}{2y^2}dz\wedge d\bar{z},
\end{equation}
while Yang \cite{Yang25} showed that
\begin{equation}
\Theta^{(2)}(\omega_0)=\sqrt{-1}\, \big(-\frac{dz\wedge d\bar{z}}{2y^2}+\frac{dw\wedge d\bar{w}}{|w|^2}+\frac{\sqrt{-1}\, dw\wedge d\bar{z}}{2yw}-\frac{\sqrt{-1}\, dz\wedge\bar{w}}{2y\bar{w}}\big)
\end{equation}
which is not non-positive.

Note that the Hermitian metric coefficients of $\omega_0$ satisfy
\[h_{w\bar{w}}=\frac{4}{|w|^2},\quad h_{w\bar{z}}=\frac{2\sqrt{-1}}{wy},\quad h_{z\bar{w}}=-\frac{2\sqrt{-1}}{y\bar{w}},\quad h_{z\bar{z}}=\frac{2}{y^2},\]
with inverse matrix
\[h^{w\bar{w}}=\frac{|w|^2}{2},\quad h^{w\bar{z}}=\frac{\sqrt{-1}\, wy}{2},\quad h^{z\bar{w}}=-\frac{\sqrt{-1}\, y\bar{w}}{2},\quad h^{z\bar{z}}=y^2.\]
Consequently,
\[
\partial^*\omega_0=-\sqrt{-1}\, \Lambda_{\omega_0}\bar{\partial}\omega_0=\sqrt{-1}\, \big(-\frac{2}{\bar{w}}d\bar{w}-\frac{\sqrt{-1}}{y}d\bar{z}\big).
\]
Moreover,
\begin{equation}\label{8.4}
\partial\partial^*\omega_0=\bar{\partial}\bar{\partial}^*\omega_0=\frac{\sqrt{-1}}{2y^2}dz\wedge d\bar{z}.
\end{equation}

By setting $t=0$ in \eqref{1.2b} and \eqref{1.2c} and applying \eqref{8.3} and \eqref{8.4}, we obtain
\begin{equation}\label{8.5}
\Theta^{(3)}(\omega_0)=\Theta^{(4)}(\omega_0)=-\frac{\sqrt{-1}}{y^2}dz\wedge d\bar{z}.
\end{equation}

Finally, taking traces in \eqref{8.3} and \eqref{8.5}, respectively, yields
\begin{equation}\label{8.6}
S_C^{(1)}(\omega_0)=-\frac{1}{2},\quad\text{and}\quad S_C^{(2)}(\omega_0)=-1.
\end{equation}

\noindent\textbf{Example 7.3.} Inoue surfaces, introduced by Inoue \cite{Ino74}, can be expressed as
$(\mathbb{H}\times\mathbb{C})/\Gamma$, where $\Gamma$ is a discrete group of
automorphisms of $\mathbb{H}\times\mathbb{C}$. They are compact, non-K\"ahler
complex surfaces of class \Rmnum{7}$_0$ with Kodaira dimension $\kappa=-\infty$,
first Betti number $b_1=1$, and second Betti number $b_2=0$.

In this example, we focus on the Inoue surface described in
\cite[Section 5]{TW13} and denote it by $S_1$.

Consider the Tricerri metric \cite{Tri82} defined by
\[
\omega_1=\frac{\sqrt{-1}}{y^2}dz\wedge d\bar{z}+\sqrt{-1}\, ydw\wedge d\bar{w},\qquad z=x+\sqrt{-1}\, y.
\]
provides a $\Gamma$-invariant Hermitian metric on $S_1$ which is Gauduchon. Tosatti and Weinkove \cite{TW13} computed that
\begin{equation}\label{8.7}
\Theta^{(1)}(\omega_1)=-\frac{\sqrt{-1}}{4y^2}dz\wedge d\bar{z}.
\end{equation}

The Hermitian metric coefficients of $\omega_1$ are
\[h_{w\bar{w}}=y,\quad h_{w\bar{z}}=h_{z\bar{w}}=0,\quad h_{z\bar{z}}=\frac{1}{y^2},\]
with inverse matrix
\[h^{w\bar{w}}=\frac{1}{y},\quad h^{w\bar{z}}=h^{z\bar{w}}=0,\quad h^{z\bar{z}}=y^2.\]
Consequently,
\[
\bar{\partial}\omega_1=-d\bar{z}\wedge dw\wedge d\bar{w},\quad\text{and}\quad
\partial^*\omega_1=-\sqrt{-1}\, \Lambda_{\omega_1}\bar{\partial}\omega_1=\frac{1}{y}d\bar{z}.
\]
Moreover,
\begin{equation}\label{8.8}
\partial\partial^*\omega_1=\bar{\partial}\bar{\partial}^*\omega_1=\frac{\sqrt{-1}}{2y^2}dz\wedge d\bar{z}.
\end{equation}

By setting $t=0$ in \eqref{1.2b} and \eqref{1.2c} and applying \eqref{8.7} and \eqref{8.8}, we obtain
\begin{equation}\label{8.9}
\Theta^{(3)}(\omega_1)=\Theta^{(4)}(\omega_1)=-\frac{3\sqrt{-1}}{4y^2}dz\wedge d\bar{z}.
\end{equation}

Finally, taking traces in \eqref{8.7} and \eqref{8.9}, respectively, yields
\begin{equation}\label{8.6}
S_C^{(1)}(\omega_1)=-\frac{1}{4},\quad\text{and}\quad S_C^{(2)}(\omega_1)=-\frac{3}{4}.
\end{equation}

\noindent\textbf{Example 7.4.} We also consider another Inoue surface following
Tosatti and Weinkove \cite[Section 6]{TW13}, which we denote by $S_2$.

For some $m\in\mathbb{R}$, the Gauduchon metric discovered by Tricerri \cite{Tri82} ($m=0$) and by Vaisman \cite{Vai87} (in general) is defined as
\[
\omega_2=\sqrt{-1}\, \big(dw-\frac{v-m\log y}{y}dz\big)\wedge\big(d\bar{w}-\frac{v-m\log y}{y}d\bar{z}\big)+\frac{\sqrt{-1}}{y^2}dz\wedge d\bar{z},
\]
where $z=x+\sqrt{-1}\, y$ and $w=u+\sqrt{-1}\, v$.

The Hermitian metric coefficients of $\omega_2$ satisfy
\[h_{w\bar{w}}=1,\quad h_{w\bar{z}}=h_{z\bar{w}}=-\frac{v-m\log y}{y},\quad h_{z\bar{z}}=\frac{1+(v-m\log y)^2}{y^2},\]
with inverse matrix
\[h^{w\bar{w}}=1+(v-m\log y)^2,\quad h^{w\bar{z}}=h^{z\bar{w}}=y(v-m\log y),\quad h^{z\bar{z}}=y^2.\]

Tosatti and Weinkove \cite{TW13} computed
\begin{equation}\label{8.10}
\Theta^{(1)}(\omega_2)=-\frac{\sqrt{-1}}{2y^2}dz\wedge d\bar{z},
\end{equation}
and
\[
\bar{\partial}\omega_2=-\frac{v-m\log y-m}{2y^2}d\bar{w}\wedge dz\wedge d\bar{z}+\frac{1}{2y}d\bar{w}\wedge dw\wedge d\bar{z}.
\]

Then we can get
\[
\partial^*\omega_2=-\sqrt{-1}\, \Lambda_{\omega_2}\bar{\partial}\omega_2=\frac{1+m(v-m\log y)}{2y}d\bar{z}-\frac{m}{2}d\bar{w},
\]
and hence
\begin{equation}\label{8.11}
\partial\partial^*\omega_2=\frac{\sqrt{-1}}{4y^2}(1+m(v-m\log y+m))dz\wedge d\bar{z}-\sqrt{-1}\, \frac{m}{4y}dw\wedge d\bar{z},
\end{equation}
Taking $t=0$ in \eqref{1.2b} and \eqref{1.2c}, and then applying \eqref{8.10} and \eqref{8.11}, respectively  gives
\begin{equation}\label{8.12}
\Theta^{(3)}(\omega_2)=-\frac{\sqrt{-1}}{4y^2}(3+m(v-m\log y+m))dz\wedge d\bar{z}+\sqrt{-1}\, \frac{m}{4y}dw\wedge d\bar{z}.
\end{equation}
Moreover, we have
\begin{equation}\label{8.13}
\Theta^{(4)}(\omega_2)=-\frac{\sqrt{-1}}{4y^2}(3+m(v-m\log y+m))dz\wedge d\bar{z}+\sqrt{-1}\, \frac{m}{4y}dz\wedge d\bar{w}.
\end{equation}

Taking traces in \eqref{8.10} and \eqref{8.12}, respectively, yields
\begin{equation}\label{8.6}
S_C^{(1)}(\omega_2)=-\frac{1}{2},\quad\text{and}\quad S_C^{(2)}(\omega_2)=-\frac{3+m^2}{4}.
\end{equation}

\section*{Acknowledgements}
The author is grateful to Kefeng Liu for his constant guidance and support, and to Jianchun Chu and Xiaokui Yang for illuminating discussions. The author also thanks Man-Chun Lee, Yong Luo and Jeffrey Streets for their interest.



\section*{Declaration of generative AI and AI-assisted technologies in the manuscript preparation process}

During the preparation of this work the author used ChatGPT in order to assist with language editing and improve the clarity of presentation. After using this tool/service, the author reviewed and edited the content as needed and takes full responsibility for the content of the published article.

\end{document}